\documentclass[11pt,reqno]{amsart}
\usepackage{graphicx, xcolor} % Required for inserting images
\usepackage[margin=2.75cm]{geometry}
\usepackage[english]{babel}

\usepackage{amssymb,amsmath,amsthm}
\usepackage{amsfonts,latexsym,amssymb,amsmath,url,amsthm, color}
\usepackage[dvipsnames]{xcolor}

\usepackage[
    colorlinks=true,
    linkcolor=Violet,
    citecolor=Violet,
    urlcolor=Violet
]{hyperref}

\usepackage{amsfonts,latexsym,url, color}
\usepackage{tikz}
\usepackage{float}
\usetikzlibrary{calc}
\newtheorem{theorem}{Theorem}[section]
\newtheorem*{theorem*}{Theorem}
\newtheorem{lemma}{Lemma}[section]

\theoremstyle{remark}
\newtheorem*{ack}{Acknowledgements}
\newtheorem*{AI}{AI statement}

\theoremstyle{definition}
\theoremstyle{remark}

\title[Primitive rank-two lattices in $\Z^4$]{Counting primitive rank-two lattices in $\Z^4$}
\author{Lena Wurzinger}
\address{Institute of Science and Technology Austria (ISTA)\\
Am Campus 1\\
3400 Klosterneuburg\\
Austria}
\email{lena.wurzinger@ist.ac.at}

\numberwithin{equation}{section}

\let\originalLambda\Lambda
\renewcommand{\Lambda}{\mathsf{\originalLambda}}

\newcommand{\Z}{\mathbb{Z}}
\newcommand{\Zprim}{\mathbb{Z}_{\mathrm{prim}}}
\newcommand{\R}{\mathbb{R}}
\newcommand{\Q}{\mathbb{Q}}

\newcommand{\C}{\mathbb{C}}

\DeclareMathOperator{\gr}{\mathbf{Gr}}

\DeclareMathOperator{\rad}{rad}


\begin{document}

\begin{abstract}
  We establish an explicit secondary term in Schmidt's asymptotic formula for the number of primitive rank-two lattices in $\Z^4$ of bounded determinant, using an exact formula for the corresponding fixed-determinant count. 
\end{abstract}

\maketitle

\setcounter{tocdepth}{1}
\tableofcontents

\thispagestyle{empty}

\section{Introduction}
      Let us write $\mathcal{L}_{k,n}$ for the set of all primitive rank-$k$ lattices in $\Z^n$  and let
      \[
      P_{k,n}(B)
      = 
      \# 
    \{ 
    \Lambda \in \mathcal{L}_{k,n} : \det(\Lambda) \leq B 
    \}, \qquad 1 \le k < n,
      \]
      denote the number of such lattices of determinant at most $B$.
      In 1968, Schmidt~\cite{schmidt_asymptotic_formulae} used the geometry of numbers to prove that 
      \begin{align}\label{eq:schmidt_asymptotic}
           P_{k,n}(B) =  c_{k,n} B^n + O\left(B^{n -\max \left\{ {1}/{k}, {1}/{(n-k)}
    \right\} }\right), 
      \end{align}
      as $B \to \infty$, 
      where
            \[
      c_{k,n} 
=
\frac{1}{n}\binom{n}{k}
\frac{V(n)}{2 \zeta(n)}
\prod_{j=2}^{k}\frac{V(n-j+1) \zeta(j)}{V(j) \zeta(n-j+1)} ,
      \]
with $ \zeta $ denoting the Riemann zeta function and $ V(d) $ the volume of the $d$-dimensional unit ball.

From an algebro-geometric perspective, the  elements of
$\mathcal{L}_{k,n}$ are naturally identified
with the rational points on the Grassmannian
 \[
    \gr (k,n) 
    = 
    \{ V: V \text{ is a $k$-dimensional subspace of } \Q^n 
    \}.
    \]
    Indeed, this identification is given by
$
    \gr(k,n)(\mathbb{Q}) \to \mathcal{L}_{k,n}, \
    V \mapsto V \cap \mathbb{Z}^n .
$
    Using the Pl\"ucker embedding 
    $ \gr(k,n) \to \mathbb{P}^{\binom{n}{k}-1}$, the Grassmannian may be realised as a smooth
projective variety.
This allows us to view results about $P_{k,n}(B)$ 
as results about counting rational points of bounded height on Fano varieties. 
Schmidt's asymptotic formula \eqref{eq:schmidt_asymptotic} counts points of bounded height on
Grassmannian varieties over $\mathbb{Q}$, verifying the Manin conjecture~\cite{manin_conjecture} for this family of varieties. As explained in ~\cite[\S~6]{peyremaninconj}, the leading constant $c_{k,n}$ agrees with Peyre's prediction. 

In the case $k=1$, estimating the error term in
\eqref{eq:schmidt_asymptotic} is the classical primitive
$n$-dimensional Gauss circle problem. 
Moreover, for any primitive lattice $\Lambda$, we have $\det(\Lambda) = \det(\Lambda^\perp)$, where $\Lambda^\perp$ denotes the \emph{orthogonal} lattice of $\Lambda$, and therefore $P_{k,n}(B) = P_{n-k,n}(B)$. The first case not covered by rank-one results is thus $(k,n)=(2,4)$. In this case, the Pl\"ucker embedding identifies $\gr(2,4)$ with the projective quadric
\begin{align}
 \label{eq:kleinquadric}
    x_1x_2 - x_3x_4 + x_5x_6 = 0,
\end{align}
known as the \emph{Klein quadric}. 
Specialising \eqref{eq:schmidt_asymptotic} gives 
\begin{equation}\label{eq:p24sum}
         P_{2,4}(B) =   \frac{15}{2 \zeta(3)}B^4 + O(B^{7/2}).
\end{equation}
We shall prove the following result in Section~\ref{section:boundedheight}.
\begin{theorem}\label{thm:boundedheight}
    As $B \to \infty$, we have 
    \begin{equation}\label{eq:sumoverfixedheight}
     P_{2,4}(B) =   \frac{15}{2 \zeta(3)}B^4 - 
     \frac{16 \log (2)}{\pi \zeta(3)} B^3
     + O_\varepsilon(B^{14/5+\varepsilon}).
    \end{equation}
\end{theorem}

Related work of Getz~\cite[Theorem~1.1]{getzsecondary} concerning smoothly weighted points on quadrics is applicable to \eqref{eq:kleinquadric}. After M\"{o}bius inversion, it yields an asymptotic of the form
\[
    \sum_{ \Lambda \in \mathcal{L}_{2,4} } f( \det(\Lambda)/B) = c_{1,f} B^{4} + c_{2,f}B^3 + O_{\varepsilon,f}(B^{2+\varepsilon})
\]
whenever $f$ satisfies suitable smoothness assumptions\footnote{In this setting, the function $\mathbf{x} \mapsto f(\|\mathbf{x}\|_2)$ needs to be smooth.} and is compactly supported. By contrast, in Theorem~\ref{thm:boundedheight} we count points with respect to the sharp cut-off function $\mathbf{1}_{[0,1]}(\det(\Lambda)/B)$.

The proof of Theorem~\ref{thm:boundedheight} relies on the identity
\[
    P_{2,4}(B)=\sum_{N\le B^2}p_{2,4}(\sqrt N),
\]
where
\[
    p_{2,4}(\sqrt N)
    =\#\{\Lambda\in\mathcal L_{2,4}:
          \det(\Lambda)=\sqrt N\}.
\]
Building on work of Aka--Einsiedler--Wieser~\cite{akaeinsiwieser}, Chinta and Pereira J\'unior~\cite[Proposition~5.4]{planesZ4} obtained an explicit formula for this fixed-determinant count. We give a short derivation of an equivalent formula
of the form 
\[
    p_{2,4}(\sqrt N)= a(N) N
\]
in Section~\ref{section:fixedproof}. Here $a(N)$ is the product of
$L(1,(\tfrac{-4N}{\cdot}))^2$ with a well-behaved arithmetic function.

Because of the structure of $a(N)$, our sum in \eqref{eq:p24sum} is amenable to methods used in the study of moments of quadratic $L$-functions. We expand and suitably truncate the Dirichlet series for the square of the $L$-value, separate the square and non-square terms and apply Poisson summation to the latter. The square terms yield the main term, while the non-oscillatory terms arising after Poisson summation produce the secondary term.

A closely related sum, analysed by similar methods, appears in the PhD thesis of Virtanen~\cite{Virtanen2008}. Virtanen obtains in \cite[Theorem~3.27]{Virtanen2008} an expansion of the form\footnote{Strictly speaking, Virtanen's \cite[Theorem~3.27]{Virtanen2008} discusses sums over all negative integers $-X\le m<0$ lying in residue classes $m\equiv 0,1 \bmod 4$. The negative odd case is discussed separately, and it is noted that the negative even case can be treated analogously.}
\[
    \sum_{N\le X}
    L\left(1,\left(\tfrac{-4N}{\cdot}\right)\right)^2
    =
    cX+P(X)\sqrt{X}
    +O_\varepsilon(X^{2/5+\varepsilon}),
\]
where $c>0$ and $P(X)$ is an explicit, potentially $X$-dependent function arising from the non-oscillatory terms after Poisson summation. In our setting, we incorporate the arithmetic weights arising from $a(N)$ into the analysis and refine the treatment of the non-oscillatory terms, obtaining an asymptotic formula for their contribution. This yields Theorem~\ref{thm:boundedheight}. 

Related work on moments of quadratic Dirichlet $L$-functions $L(s,\chi)$ often concerns the central point $s=1/2$. In particular, the treatment of the non-square terms
in Virtanen's argument builds on the work of Conrey, Farmer and
Soundararajan~\cite{conreyfarmersoundararajan2000} on double sums of quadratic
characters, which arose in connection with Soundararajan's work on the
non-vanishing of quadratic Dirichlet $L$-functions at the central
point~\cite{soundnonvanishing}. Results more directly related to the value $s=1$ include work of Wolke~\cite{WOLKE1972523} and Jutila~\cite{JUTILA1973203} on second moments of class numbers associated to negative discriminants. Granville and Soundararajan~\cite{granvillesoundararajan2003} studied the distribution of
$
L\left(1,\left(\tfrac{d}{\cdot}\right)\right)
$
as $d$ varies over fundamental discriminants, including its moments.

 \subsection*{Definitions and notation} 
 Let $w_1, \ldots, w_k \in \Z^n$ be linearly independent vectors. We call the set of all linear integer combinations of these vectors, 
    $
    \Lambda =  \mathrm{span}_{\Z}(w_1, \ldots, w_k),
    $
    an \emph{integral lattice of rank $k$}. 
    An integral lattice is called \emph{primitive}  if 
    $
     \mathrm{span}_{\Q}(\Lambda) \cap \Z^n = \Lambda. 
    $
    The \emph{determinant} of $\Lambda$ is 
    \begin{equation*}
    % \label{eq:defdet}
        \det(\Lambda) =
   \sqrt{\det (M_\Lambda M_\Lambda^t)},
    \end{equation*}
        where
     $M_\Lambda$ is a matrix whose rows form a basis of $\Lambda$
     and
    $M_\Lambda^t$ denotes its transpose. 
     One may check that the determinant is independent of the choice of basis.

For any positive integer $N$, we set
\begin{equation*}
% \label{eq:def_lfnc}
    L\left(s,\left(\tfrac{-4N}{\cdot}\right)\right)
    =
    \sum_{n=1}^{\infty}
    \frac{1}{n^s}
    \left(\frac{-4N}{n}\right) = \prod_{p}
    \left(1-\left(\frac{-4N}{p}\right)p^{-s}\right)^{-1},
\end{equation*}
where $\left(\frac{\cdot}{\cdot}\right)$ denotes the Kronecker symbol. 
Both the series and the Euler product converge absolutely for $\Re(s)>1$;
at $s=1$, they converge conditionally to the same value. 
We write $\rad(n)=\prod_{p \mid n} p $ for the \emph{radical} of a positive integer $n$.

We use $O(\cdot)$ for Big--$O$ notation and $\ll,\gg$ in the usual
Vinogradov sense. Dependencies of implied constants are indicated by subscripts, as in
$O_\alpha(\cdot)$ or $\ll_\alpha$. Finally, $\varepsilon>0$ always denotes
an arbitrarily small positive constant, which may vary from line to line.

\begin{ack}
    The author is very thankful to her advisor Tim Browning for suggesting the project and his guidance. She would like to thank Damaris Schindler and Alistair Severn for interesting discussions, and she would like to thank Nikita P.~Kalinin for his assistance with AI tools.
\end{ack}  

\begin{AI}
    Generative AI tools were used to create an outline for the proof of Theorem~\ref{thm:boundedheight}. Specifically, GPT-5.6 Sol and Claude Opus 5 were used in order to improve upon an argument in an earlier draft, which gave an error term of $O_\varepsilon (B^{3+\varepsilon})$.  GPT-6.2  was used for proofreading of the entire manuscript. All arguments were independently checked by the author, who takes full responsibility for the contents of the paper. The final text was written by the author.
\end{AI}

\section{Primitive lattices of fixed determinant} \label{section:fixedproof}

In this section, we give our own proof of the following identity.

    \begin{theorem}[{Chinta and Pereira J\'unior \cite[Proposition~5.4]{planesZ4}}] \label{thm:2-4}
    Let $N$ be a positive integer.
        Then $p_{2,4}(\sqrt{N}) = a(N) N$, where
    \[
    a(N) = \frac{b_N}{\pi^2} 
    L \left(1,  \left(\tfrac{-4 N }{ \cdot } \right) \right)^2  
          \prod_{\substack{p^2 \mid N, \\ p\neq 2}} (1+g_{N}(p)),
\]
% with the product going over all odd primes such that $p^2 \mid N$,
with 
$$g_N(p) =  
        {2}{p^{-\lfloor\frac{\nu_p(N)}{2}\rfloor}} 
        \left( 
        \frac{ 
         p^{\lfloor\frac{\nu_p(N)}{2}\rfloor} -p
        }{
        p-1
        } +
      \left( 
            1 - p^{-1}
            \Bigg(\frac{- N / p^{2\lfloor\frac{\nu_p(N)}{2}\rfloor}}{ p } 
        \right)
         \right)^{-1}
         \Bigg),
$$
and 
\[
b_N = \begin{cases}
    0 &\text{if } N \equiv 0, 7, 12, 15  \bmod 16, \\
    128 &\text{if } N \equiv 3, 11  \bmod 16, \\
    96 &\text{if } N \equiv 1,2,5,6,9, 10, 13, 14  \bmod 16, \\
    48 &\text{if } N \equiv  4,8  \bmod 16.
\end{cases}
\]
\end{theorem}

Chinta and Pereira J\'unior used their formula to obtain
explicit expressions for the Dirichlet series counting 
primitive rank-two lattices in $\Z^4$ and to study its 
relation to Eisenstein series on $\mathrm{GL}_4$. 
We begin by proving the following lemma.

\begin{lemma}\label{lem:prim}
    We have  
    \[
p_{2,4}(\sqrt{N})
=
\frac{1}{2}\# \mathcal{A}_{N}^\mathrm{prim} ,
\]
where 
\[
\mathcal{A}_{N}^\mathrm{prim} =
\left\{
(x_1,\ldots,x_6)\in \Zprim^6
:~
x_1x_2 - x_3x_4 + x_5x_6 = 0,~
x_1^2+x_2^2+x_3^2+x_4^2+x_5^2+x_6^2 = N
\right\}.
\]
\end{lemma}

\begin{proof}
As discussed in the introduction, the Pl\"ucker embedding
$\iota \colon \gr(2,4) \hookrightarrow \mathbb{P}^5$
identifies $\gr(2,4)$ with the projective quadric
\[
    z_{1,2}z_{3,4} - z_{1,3}z_{2,4} + z_{1,4}z_{2,3} = 0.
\]
 More explicitly, for any two-dimensional vector space $V = \operatorname{span}_{\Q}(\mathbf{v},\mathbf{w})$ with basis vectors $\mathbf{v}=(v_1,v_2,v_3,v_4)$ and $\mathbf{w}=(w_1,w_2,w_3,w_4)$, its image under the Pl\"ucker embedding $\iota$, independently of the choice of basis, is the collection of minors 
\[
\left(
\begin{vmatrix}
v_1 & v_2 \\
w_1 & w_2
\end{vmatrix}
:
\begin{vmatrix}
v_1 & v_3 \\
w_1 & w_3
\end{vmatrix}
:
\begin{vmatrix}
v_1 & v_4 \\
w_1 & w_4
\end{vmatrix}
:
\begin{vmatrix}
v_2 & v_3 \\
w_2 & w_3
\end{vmatrix}
:
\begin{vmatrix}
v_2 & v_4 \\
w_2 & w_4
\end{vmatrix}
:
\begin{vmatrix}
v_3 & v_4 \\
w_3 & w_4
\end{vmatrix}
\right).
\]

Moreover, there is a natural bijection
\begin{align*}
    \gr(2,4)(\mathbb{Q}) &\longrightarrow \mathcal{L}_{2,4}, \\
    V &\longmapsto V \cap \mathbb{Z}^4
\end{align*} 
between primitive rank-two sublattices of $\mathbb{Z}^4$ and two-dimensional $\mathbb{Q}$-vector subspaces of $\mathbb{Q}^4$.
If the Pl\"ucker coordinates $
    \iota(V) = (y_{1,2}:y_{1,3}:y_{1,4}:y_{2,3}:y_{2,4}:y_{3,4})
$ are given as integers $y_{i,j}$ with 
$
     (y_{1,2},\ldots,y_{3,4})=1,
$
 then it may be shown that 
\begin{equation} \label{eq:detsumofsquares}
    \det(V \cap \mathbb{Z}^4)^2
    =
    y_{1,2}^2+y_{1,3}^2+y_{1,4}^2+y_{2,3}^2+y_{2,4}^2+y_{3,4}^2.
\end{equation}
 Every lattice thus corresponds to two elements of $\mathcal{A}_N^{\mathrm{prim}}$, namely 
$
    \pm (y_{1,2},\ldots,y_{3,4}).
$
The equality in \eqref{eq:detsumofsquares} between the determinant and the Euclidean norm of the Pl\"ucker coordinates holds more generally for Grassmannians.

For the proof,  if $\mathbf{v}$ and $\mathbf{w}$ span a primitive lattice, then the minors 
\[
\left(
\begin{vmatrix}
v_1 & v_2 \\
w_1 & w_2
\end{vmatrix}
,
\begin{vmatrix}
v_1 & v_3 \\
w_1 & w_3
\end{vmatrix}
,
\begin{vmatrix}
v_1 & v_4 \\
w_1 & w_4
\end{vmatrix}
,
\begin{vmatrix}
v_2 & v_3 \\
w_2 & w_3
\end{vmatrix}
,
\begin{vmatrix}
v_2 & v_4 \\
w_2 & w_4
\end{vmatrix}
,
\begin{vmatrix}
v_3 & v_4 \\
w_3 & w_4
\end{vmatrix}
\right)
\]
form a primitive vector $(y_{1,2},\ldots,y_{3,4})\in\Zprim^6$.  Then \eqref{eq:detsumofsquares} follows by the Cauchy--Binet formula.

In order to see that $ (y_{1,2},\ldots,y_{3,4})=1$, assume that there is a prime $p$ such that $p$ divides all $y_{i,j}$. As a basis vector of a primitive lattice is primitive, we may assume without loss of generality that $p\nmid v_1$. Moreover, after a potential change of basis by adding a suitable multiple of $\mathbf{v}$ to $\mathbf{w}$, we may assume that $p \mid w_1$. But then for every $j\ge 2$ we have
\[
    p \mid w_1 v_j - v_1 w_j.
\]
As $p \mid w_1$ and $p \nmid v_1$, it follows that $p \mid w_j$ and so each coordinate of $\mathbf{w}$ is divisible by $p$, which contradicts the fact that $\mathbf{w}$ is a basis vector of a primitive lattice and thus primitive itself.
\end{proof}

\begin{lemma} \label{lem:all_solutions_correspondence}
Let $N$ be a positive integer and let 
\begin{align*}
    \mathcal{A}_N &= 
\left\{
(x_1,\ldots,x_6)\in \Z^6
: 
x_1x_2 - x_3x_4 + x_5x_6 = 0,~
x_1^2+x_2^2+x_3^2+x_4^2+x_5^2+x_6^2 = N
\right\}, \\
 \mathcal{B}_N &=
\left\{
(y_1,\ldots,y_6)\in \Z^6
:~
\begin{aligned}
y_1^2 + y_4^2 + y_5^2 = y_2^2 + y_3^2 + y_6^2 = N
\end{aligned}
\right\}.
\end{align*}
Then 
\[
\# \mathcal{A}_N = C_N \# \mathcal{B}_N,
\]
where
    \begin{equation} \label{eq:defCm}
    C_N= \begin{cases}
        1  & \text{if } N \equiv 0,3 \bmod 4, \\
        \frac{1}{3}   & \text{if }   N \equiv 1,2 \bmod 4. 
    \end{cases}
    \end{equation}
\end{lemma}
\begin{proof}
First, let us choose $(x_1, \ldots, x_6) \in  \mathcal{A}_N$.
We apply the following change of coordinates: 
\begin{equation*}
    \begin{matrix}
        y_1 = x_1 + x_2, \qquad y_3 = x_3 + x_4, \qquad y_5 = x_5 + x_6, \\
         y_2 = x_1 - x_2, \qquad 
         y_4 = x_3 - x_4, \qquad
        y_6 = x_5 - x_6. 
\end{matrix}
\end{equation*}
As $x_1x_2 - x_3x_4 + x_5x_6 = 0$, we have
\begin{align*}
    y_1^2 + y_4^2 + y_5^2 = 
    x_1^2 + x_2^2 + x_3^2 + x_4^2 + x_5^2 + x_6^2 
    = N
\end{align*} and 
\begin{align*}
    y_2^2 + y_3^2 + y_6^2 =
    x_1^2 + x_2^2 + x_3^2 + x_4^2 + x_5^2 + x_6^2 
    = N;
\end{align*}
 hence it follows that $(y_1, y_2, y_3, y_4, y_5, y_6) \in  \mathcal{B}_N$.

Conversely, given $(y_1, \ldots, y_6) \in  \mathcal{B}_N$, 
we can find $(x_1, \ldots, x_6) \in \mathbb{Q}^6$ with   
\begin{align*}
    x_1x_2 - x_3x_4 + x_5x_6 = 0 
    \quad 
    \text{ and }
    \quad
x_1^2+x_2^2+x_3^2+x_4^2+x_5^2+x_6^2 = N
\end{align*}
 via $x_k = (y_k + y_{k+1})/2$, $x_{k+1} = (y_k - y_{k+1})/2$ for $k = 1,3,5$. 
It is clear that all $x_i$ are integers, i.e. that $(x_1,x_2,x_3,x_4,x_5,x_6)\in \mathcal{A}_N$, if and only if the two integers in each pair $y_1$ and $y_2$, $y_3$ and $y_4$ as well as $y_5$ and $y_6$ have the same 
parity.

If $N \equiv 0,3 \bmod 4$, then considering the congruences 
\begin{equation*}
    y_1^2 + y_4^2 + y_5^2 \equiv N \bmod 4, \qquad
    y_2^2 + y_3^2 + y_6^2 \equiv N \bmod 4
\end{equation*}
implies that all $y_i$ are of the same parity; indeed, all $y_i$ are even if $N \equiv 0 \bmod 4$ and all $y_i$ are odd if $N \equiv 3 \bmod 4$. Therefore the change of variables outlined above yields a bijection between the sets $\mathcal{A}_N$ and $\mathcal{B}_N$ and so 
$
    \# \mathcal{A}_N = \# \mathcal{B}_N.
$

If $N \equiv 1,2 \bmod 4$, then considering the congruences 
\begin{equation*}
    y_1^2 + y_4^2 + y_5^2 \equiv N \bmod 4, \qquad
    y_2^2 + y_3^2 + y_6^2 \equiv N \bmod 4
\end{equation*}
implies that, per equation, exactly two of the $y_i$ are of the same parity; indeed, exactly one $y_i$ per equation is odd if $N \equiv 1 \bmod 4$ and exactly one $y_i$ per equation is even if $N \equiv 2 \bmod 4$. By symmetry, for exactly one third of all pairs of triples $(y_1, y_4, y_5)$ and $(y_2,y_3,y_6)$ the uniquely (in terms of parity) distinguished entry occurs in the same position in both triples. It follows that
$
    \# \mathcal{A}_N = \frac{1}{3} \# \mathcal{B}_N.
$
\end{proof}

\begin{lemma}[Gauss, see e.g.~{\cite[Theorem B]{bateman1951representations}}] \label{thm:bateman}
    Let $N$ be a positive integer, let 
    $$
    R_3(N) = 
    \# 
    \{
    (x_1, x_2, x_3) \in \Zprim^3 :~ x_1^2 + x_2^2 + x_3^2 = N
    \}
    $$ denote the number of primitive representations of $N$ as a sum of three squares and let 
    $$
    r_3(N) = \# 
    \{
    (x_1, x_2, x_3) \in \Z^3 :~ x_1^2 + x_2^2 + x_3^2 = N
    \}
    $$
    denote the number of all representations of $N$ as a sum of three squares.

    Then we have
    \[
    R_3(N) = \frac{G_N}{\pi}\sqrt{N} L\left(1, \left(\tfrac{-4N}{\cdot}\right)\right) 
    \]
    and 
    \[
        r_3(N) = \sum_{d^2 \mid N}
        R_3\left({N}/{d^2} \right),
    \]
     where   
    \[
    G_N = \begin{cases}
        0  & \text{if } N \equiv 0,4,7 \bmod 8, \\
        16 & \text{if }   N \equiv 3 \bmod 8, \\
        24 & \text{if }  N \equiv 1,2,5,6 \bmod 8  .
    \end{cases}
    \]
\end{lemma}

\begin{proof}[Proof of Theorem~\ref{thm:2-4}]
Define $\mathcal{A}_N^{\mathrm{prim}}$ and $\mathcal{A}_N$ as in Lemma~\ref{lem:prim} and Lemma~\ref{lem:all_solutions_correspondence}. By Lemma~\ref{lem:prim} and M\"{o}bius inversion, we have  
\begin{align*}
    2 p_{2,4}(\sqrt{N}) 
   = \sum_{d^2 \mid N} \mu (d)\# \mathcal{A}_{N/d^2}.
\end{align*} Lemma~\ref{lem:all_solutions_correspondence} yields   
$
   \# \mathcal{A}_{m} = C_m r_3(m)^2,
$
for all $m\ge 1$, where $C_m$ is defined in \eqref{eq:defCm}.
Lemma~\ref{thm:bateman} gives
\begin{align}
    2p_{2,4}(\sqrt{N}) 
    &= \! \sum_{d^2 \mid N} \mu (d) C_{{N}/{d^2}}
         r_3(N/d^2)^2  \label{eq:asbefore}  \\
         &= \! \sum_{d^2 \mid N} \mu (d)C_{\frac{N}{d^2}} 
         \bigg( \sum_{t^2 \mid  \frac{N}{d^2}} 
        R_3\left({N}/{(dt)^2} \right) \bigg)^2  \nonumber \\
        &= \! \sum_{d^2 \mid N} \mu (d)C_{{N}/{d^2}}
         \sum_{t_1^2, t_2^2 \mid  \frac{N}{d^2}} 
        R_3\left({N}/{(d t_1)^2} \right)
        R_3\left({N}/{(dt_2)^2} \right)  \nonumber \\
        &= \! \sum_{d^2 \mid N} \mu (d)C_{{N}/{d^2}}
         \!\sum_{t_1^2, t_2^2 \mid  \frac{N}{d^2}} 
         \frac{N}{\pi^2t_1 t_2 d^2}
         \prod_{i=1}^2
         \left(
         G_{{N}/{(dt_i)^2}} L\left(1,\left(\tfrac{{-4N}/{(dt_i)^2}}{\cdot}\right) \right)  
         \right). \nonumber  
         \end{align}

          Let us first assume that $N$ is of the form $N=2^k N'$ with $k \in \{0,1\}$ and $N'$ odd.
          In this case, any $d$ with $d^2|N$ is odd, hence $N/d^2 \equiv N \mod 4$;  we thus get 
        $
       C_{{N}/{d^2}}= C_N
        $
        as well as 
        $
        G_{{N}/{d^2}} = G_N.
        $ Moreover, we observe that 
        \[
        L\left(
        1, \left(
        \tfrac{-4N/d^2}{\cdot}\right)
        \right) = 
         \prod_{p \mid N}
        \left( 
            1 - \left(\tfrac{- 4N / d^2 }{ p } 
        \right)p^{-1}
         \right)^{-1}
        L\left(
        1, \left(
        \tfrac{-4N}{\cdot}\right)
        \right).
        \]
        Together, this yields 
\begin{align} \label{eq:fullsum}
        2p_{2,4}(\sqrt{N})&=N \frac{C_N G_N^2}{\pi^2} L(1, \left( \tfrac{-4N}{\cdot}
        \right))^2 
         \sum_{d^2 \mid N}
         \sum_{t_1^2, t_2^2 \mid  \frac{N}{d^2}}
         \frac{\mu (d)}{t_1 t_2 d^2} 
         \prod_{i=1}^2
          \prod_{p \mid N}
        \left( 
            1 - \left(\tfrac{- 4N / (d^2 t_i^2)}{ p } 
        \right) p^{-1}
         \right)^{-1}. 
\end{align}

 We let 
 \[
f(d,t_1, t_2) := \frac{\mu (d)}{t_1 t_2 d^2} 
         \prod_{i=1}^2
          \prod_{p \mid N}
        \left( 
            1 - \left(\tfrac{-4N / (d^2 t_i^2)}{ p } 
        \right) p^{-1}
         \right)^{-1},
\]
so that it suffices to compute 
\begin{equation}\label{eq:sum_with_d1}
     \sum_{d^2 \mid N}
         \sum_{t_1^2, t_2^2 \mid  \frac{N}{d^2}} f(d,t_1, t_2).
\end{equation}
The function $f$ is jointly multiplicative; indeed, for all summation variables $(t_1d)^2, (t_2 d)^2 \mid N$ and  $(s_1e)^2, (s_2 e)^2 | N$ 
such that  $(t_1 t_2 d, s_1 s_2 e) =1$ we have 
\[
    f(de, t_1 s_1, t_2 s_2) = f(d,t_1,t_2) f(e,s_1,s_2). 
\]
Therefore 
\[
    \sum_{\substack{t_1^2d^2 \mid N, \\ t_2^2 d^2 \mid N}}
          f(d,t_1, t_2) 
         = \prod_{p^2 \mid N} \sum_{\substack{k_1 + \ell \le \nu_p(N)/2, \\ k_2 + \ell \le \nu_p(N)/2 }} f(p^\ell, p^{k_1}, p^{k_2}).
\]

Now, if $\ell \ge 2$, then
\[
    f(p^\ell,p^{k_1}, p^{k_2}) =0
\]
for any $k_1,k_2$ since $\mu(p^\ell)=0$. Moreover, for any tuples $0 \le k_1, k_2\le \nu_p(N)/2 -1$ we have 
\[
    f(p, p^{k_1},p^{k_2}) = -f(1, p^{k_1+1},p^{k_2+1}).
\]
This means that all summands corresponding to indices of the form $(1,k_1,k_2)$ are cancelled out by summands corresponding to indices $(0,k_1 + 1, k_2 +1)$. 
The only remaining summands correspond to indices $(0, k_1, 0)$ or $(0, 0,k_2)$ with $0 \le k_1, k_2 \le \nu_p(N)/2$, with $(0,0,0)$ counted only once. Hence
\[
     \prod_{p^2 \mid N} \sum_{\substack{k_1 + \ell \le \nu_p(N)/2, \\ k_2 + \ell \le \nu_p(N)/2 }} f(p^\ell, p^{k_1}, p^{k_2}) = \prod_{p^2 \mid N} \left(
     1 + 2 \sum_{1 \le k \le \nu_p(N)/2} \frac{1}{p^k}  \left( 
            1 - \left(\tfrac{-4N / p^{2k}}{ p } 
        \right) p^{-1}
         \right)^{-1}
     \right)
\]
Since $\left(\tfrac{-4N / p^{2k}}{ p } 
        \right) = 0$ for $k < \nu_p(N)/2$, we may simplify this as
\begin{align}\label{eq:sumrepresentationofgN}
    \prod_{p^2 \mid N} \left(
     1 + 2 \sum_{1 \le k \le \nu_p(N)/2} \frac{1}{p^k}  \left( 
            1 - \left(\tfrac{-4N / p^{2k}}{ p } 
        \right) p^{-1}
         \right)^{-1}
     \right) = \prod_{p^2 \mid N}
   (1+ g_N(p)),
\end{align}
where 
\[
    g_N(p) 
    =
        2{p^{-\lfloor\frac{\nu_p(N)}{2}\rfloor}} \left( 
        \frac{ 
         p^{\lfloor\frac{\nu_p(N)}{2}\rfloor} -p
        }{
        p-1
        } +
      \Bigg( 
            1 - p^{-1}\left(\tfrac{- N / p^{2\lfloor\frac{\nu_p(N)}{2}\rfloor}}{ p } 
         \right) 
          \Bigg)^{-1}
         \right).
\]
Substituting back into \eqref{eq:fullsum} yields
\[
 2p_{2,4}(\sqrt{N})=N \frac{C_N G_N^2}{\pi^2} L(1, \left( \tfrac{-4N}{\cdot}
        \right))^2 
        \prod_{p^2 \mid N} (1+g_N(p)).
\]
Using the definition of $C_N$ from \eqref{eq:defCm} 
% and the definition of $G_N$ from Lemma~\ref{thm:bateman} 
gives 
\[
   p_{2,4}(\sqrt{N})= \frac{b_N}{\pi^2}N  L(1, \left( \tfrac{-4N}{\cdot}
        \right))^2 
        \prod_{p^2 \mid N} (1+g_N(p)),
\]
where 
\[
b_N = \frac{C_N G_N^2}{2} = \begin{cases}
    0 &\text{if } N \equiv  7, 12, 15  \bmod 16, \\
    128 &\text{if } N \equiv 3, 11  \bmod 16, \\
    96 &\text{if } N \equiv 1,2,5,6,9, 10, 13, 14  \bmod 16, 
\end{cases}
\]
and thus concludes the proof for $N$ with $\nu_2(N)\le 1$. 

\medskip 

Now let us consider the case where $4 \mid N$, writing $N = 4^k M$ with $4 \nmid M$. 
By \eqref{eq:asbefore}, we have  
\begin{equation}\label{eq:p24generalN}
     2 p_{2,4}(\sqrt{N}) = \sum_{d^2 \mid N} \mu (d)C_{{N}/{d^2}}
         r_3(N/d^2)^2.
\end{equation}
Noting that, for every odd $d$ with $d^2\mid N$ (as $R_3(K)=0$ if $4 \mid K$ by Lemma~\ref{thm:bateman}),
$
    r_3(N/d^2) 
    = r_3(M/d^2)
$
and that $\mu(4d')=0$ for any integer $d'$, 
we may rewrite \eqref{eq:p24generalN} as
\begin{align*}
       p_{2,4}(\sqrt{N}) 
       &= \frac{1}{2}\sum_{\substack{d^2 \mid N,\\ d \text{ odd}}} \big( \mu (d)C_{{N}/{d^2}}
         +
         \mu (2d)C_{{N}/{(4d^2)}}
         \big) r_3(M/d^2)^2.
\end{align*}
Since $\mu(2d) = - \mu(d)$ and $C_m = C_{m/d^2}$ for all odd $d^2 \mid m$ by the definition of $C_m$ in \eqref{eq:defCm}, we get 
\begin{align*}
    p_{2,4}(\sqrt{N})
    &=
    \frac{C_N-C_{N/4}}{2}
    \sum_{d^2 \mid M} \mu(d) r_3(M/d^2)^2 \\
    &= \frac{C_N-C_{N/4}}{C_M}
    p_{2,4}(\sqrt{M}).
\end{align*}
By the definition of $C_m$ in \eqref{eq:defCm} and the definition of $G_m$ in Lemma~\ref{thm:bateman}, 
\[
    C_{{N}}-C_{\frac{N}{4}} = 1 - C_{\frac{N}{4}} =
    \begin{cases}
        0 & \text{ if } N \equiv 0, 12 \bmod 16, \\ 
        \frac{2}{3} & \text{ if } N \equiv 4, 8 \bmod 16.
    \end{cases}
\]
This concludes the proof for all $N$.
 \end{proof}

\section{Bounded-height asymptotics} \label{section:boundedheight}

Recall from Theorem~\ref{thm:2-4} that 
\[
    a(N) = \frac{b_N}{\pi^2} L(1,\left( \tfrac{-4N}{\cdot} \right))^2 \prod_{\substack{p^2 \mid N, \\ p \neq 2}} 
    \left( 1+g_N(p)
    \right),
\]
where $b_N$ depends on $N \bmod 16$ and $g_N(p)$ 
only depends on $\nu_p(N)$ and $\left( \tfrac{-N/p^{2 \lfloor \nu_p(N)/2 \rfloor}}{p} \right)$.
We will prove the following main lemma. 

\begin{lemma}\label{lem:mainlemmaboundedheight} As $X \to \infty$,  
\begin{align*}
    \sum_{N \le X} a(N) = \frac{15}{\zeta(3)} X -
    \frac{48 \log (2)}{\pi \zeta(3)}\sqrt{X} + O_\varepsilon(X^{2/5 + \varepsilon}). 
\end{align*}  
\end{lemma}

For the proof of Lemma~\ref{lem:mainlemmaboundedheight} we will use five estimates, established in the sections below. 
\begin{enumerate}
    \item In Section~\ref{section:trncating_sq_contr}, we truncate the Dirichlet series   
    \[L(1,\left( \tfrac{-4N}{\cdot} \right))^2 = \sum_{k\le Y} \left( \tfrac{-4N}{k} \right) \tfrac{d(k)}{k} 
    + \sum_{k > Y} \left( \tfrac{-4N}{k} \right) \tfrac{d(k)}{k},
    \]
    introducing some truncation parameter $Y \ge 1$. Bounding the contribution from $k > Y$ will produce
    \[
         \sum_{N \le X} a(N) = \sum_{N \le X} \frac{b_N}{\pi^2}
         \prod_{\substack{p^2 \mid N, \\ p \neq 2}} 
    \left( 1+g_N(p)
    \right)
         \sum_{k\le Y} \left( \tfrac{-4N}{k} \right) \tfrac{d(k)}{k}   + O_\varepsilon(X^{1+\varepsilon} Y^{-1/2+\varepsilon}).
    \]
    \item In Section~\ref{section:asympt_sq_contribution}, we evaluate the contribution from all $k \le Y$ that are squares, giving 
    \[
     \sum_{N \le X}
     \frac{b_N}{\pi^2} 
     \prod_{\substack{p^2 \mid N, \\ p \neq 2}} 
    \left( 1+g_N(p)
    \right)
     \sum_{
     \substack{k \le Y, \\ k = \square
     }
     }
     \frac{d(k)}{k}
    \left( \frac{-4N}{k} \right)
      = \frac{15}{\zeta(3)} X + O_\varepsilon(X^\varepsilon + XY^{-1/2+\varepsilon}).
    \]
    For the remainder of the proof, we thus only need to deal with the contribution from $k \neq \square$.
    \item In Section~\ref{section:applying_poisson}, we replace the sharp cut-off $\mathbf{1}_{[N \leq X]}$ by a smooth weight $w(X^{-1}N)$, with the  derivatives of $w$ depending on a new parameter $Z \in (2, X]$. The error introduced by smoothing will be of order 
    \[
        \ll_\varepsilon 
         (XYZ)^\varepsilon X^{1/2} Z^{-1/4}.
    \]
    Moreover, we will apply Poisson summation to our smoothed sum.
    \item In Section~\ref{section:asymptotics_after_poisson}, we evaluate the contribution from the terms that are non-oscillatory after Poisson summation, which gives
    \[
        \frac{-48 \log(2)}{\pi \zeta(3)} \sqrt{X}
      + O_\varepsilon \big( 
      (XYZ)^\varepsilon (X Y^{-1/2} + Z + (X/Z)^{1/2})
      \big).
    \]
    \item In Section~\ref{section:remaining_error_term_intheend}, we bound the remaining contribution from the oscillatory terms. This gives an additional error term of order 
    \[
    \ll_\varepsilon (XYZ)^\varepsilon Z.
    \]
\end{enumerate}
Combining these results with 
\begin{align}
     \label{eq:choiceofyz}
    Y = X^{6/5} \qquad \text{and} \qquad Z = X^{2/5}
\end{align}
yields the asymptotic claimed in Lemma~\ref{lem:mainlemmaboundedheight}, i.e.
\begin{align*}
    \sum_{N \le X} a(N) = \frac{15}{\zeta(3)} X -
    \frac{48 \log (2)}{\pi \zeta(3)}\sqrt{X} + O_\varepsilon(X^{2/5 + \varepsilon}).
\end{align*} 

\begin{proof}[Proof of Theorem~\ref{thm:boundedheight}]
    Once established, Lemma~\ref{lem:mainlemmaboundedheight} will immediately imply Theorem~\ref{thm:boundedheight} by partial summation. Indeed, we have 
\begin{align*}
        P_{2,4}(B) 
        = \sum_{N \le B^2} a(N) N 
        &= B^2 \sum_{N \le B^2} a(N) - \int_{1}^{B^2} \sum_{N \le t} a(N) \mathrm{dt} \\
        &= B^2 
        \left( 
         \frac{15}{\zeta(3)} B^2 -
    \frac{48 \log (2)}{\pi \zeta(3)}B + O_\varepsilon(B^{4/5 + \varepsilon})
        \right) \\
        &- \int_1^{B^2} \left( 
        \frac{15}{\zeta(3)} t -
    \frac{48 \log (2)}{\pi \zeta(3)} \sqrt{t} + O_\varepsilon(t^{2/5 + \varepsilon})
        \right) \mathrm{d}t \\
        & = \frac{15}{2\zeta(3)}B^4 - \frac{16 \log(2)}{\pi \zeta(3)}B^3 + O_\varepsilon(B^{14/5+\varepsilon}),
\end{align*}
which gives the claimed asymptotic. \end{proof}

\subsection{Truncating $L(1,\left( \tfrac{-4N}{\cdot} \right))^2$}\label{section:trncating_sq_contr}
In this section, we prove the following lemma.

\begin{lemma}\label{lem:truncation}
Let $Y \ge 1$. Then we have 
\begin{align}\label{eq:asumdecompwithY}
    \sum_{N \le X} a(N) = \sum_{k \le Y} \sum_{N \le X} \frac{b_N}{\pi^2} G(N)
    \left( \tfrac{-4N}{k} \right) \tfrac{d(k)}{k}
    + O_\varepsilon (X^{1+\varepsilon}Y^{-1/2+\varepsilon}),
\end{align}
where we set
\begin{align}\label{eq:defnofG}
    G(N) = \prod_{\substack{p^2 \mid N, \\ p\neq 2}} (1+g_{N}(p)).
\end{align}
\end{lemma}

\begin{proof}
    We truncate  
    \begin{align*}
        L(1,\left( \tfrac{-4N}{\cdot} \right))^2 
        =\sum_{m \le Y} \left( \tfrac{-4N}{m} \right) \tfrac{d(m)}{m}  + \sum_{m > Y} \left( \tfrac{-4N}{m} \right) \tfrac{d(m)}{m}.
    \end{align*}
    As we have 
    \begin{align*}
        \sum_{N \le X} a(N) &= \sum_{N \le X} \frac{b_N}{\pi^2} L(1,\left( \tfrac{-4N}{\cdot} \right))^2 G(N) \\
        & = 
        \sum_{N \le X} \sum_{k \le Y}
        \frac{b_N}{\pi^2} G(N) \left( \tfrac{-4N}{k} \right) \tfrac{d(k)}{k}  + 
        \sum_{N \le X} \sum_{k >  Y}
        \frac{b_N}{\pi^2} G(N) \left( \tfrac{-4N}{k} \right) \tfrac{d(k)}{k},
    \end{align*}
    the lemma statement follows immediately once we establish
    \[
     \sum_{N \le X} \sum_{k >  Y}
        \frac{b_N}{\pi^2} G(N) \left( \tfrac{-4N}{k} \right) \tfrac{d(k)}{k} = O(X^{1+\varepsilon}Y^{-1/2+\varepsilon}). 
    \]

    We decompose the inner sum dyadically. As $g_N(p) \ll p^{-1}$, we have 
    \[
    G(N) = \prod_{p^2 | N,\, p \neq 2}(1+g_N(p)) \ll N^\varepsilon
    \]
    by the divisor bound, 
    and so \cite[Lemma~3.2]{Virtanen2008} gives 
    \begin{align*}
    \sum_{N \le X} \frac{b_N}{\pi^2} G(N) \bigg| \sum_{T \le k < 2T}
    \left( \tfrac{-4N}{k} \right) \tfrac{d(k)}{k}\bigg| &\ll_\varepsilon
     X^\varepsilon \sum_{N \le X}  \bigg| \sum_{T \le k < 2T}
    \left( \tfrac{-4N}{k} \right) \tfrac{d(k)}{k}\bigg|
    \\
    &\ll_\varepsilon X^{1+2\varepsilon} T^{-1/2+\varepsilon}.
    \end{align*}
    This shows 
    \begin{align*}
         \sum_{N \le X} \sum_{k >  Y}
        \frac{b_N}{\pi^2} G(N) \left( \tfrac{-4N}{k} \right) \tfrac{d(k)}{k}  &= \sum_{j \ge 0}  \sum_{N \le X} \frac{b_N}{\pi^2} G(N)  
        \sum_{2^j Y < k \le 2^{j+1}Y}
         \left( \tfrac{-4N}{k} \right) \tfrac{d(k)}{k} \\
         &\ll_\varepsilon X^{1+2\varepsilon} Y^{-1/2+\varepsilon} \sum_{j \ge 0}  2^{-j/3} \\
         &\ll_\varepsilon X^{1+2\varepsilon} Y^{-1/2+\varepsilon}. 
    \end{align*}
    Here we took $\varepsilon<1/6$, so that we have
$2^{-j/2+j\varepsilon}\le 2^{-j/3}$ in the dyadic sum.
\end{proof}

\subsection{Asymptotics for the square contribution} \label{section:asympt_sq_contribution}
We now evaluate the contribution in \eqref{eq:asumdecompwithY} from all $k$ that are a perfect square.

\begin{lemma}\label{lem:squarecontribution}
    We have 
    \begin{align}\label{eq:squaresum}
     \sum_{
     \substack{k \le Y, \\ k = \square
     }
     }
     \sum_{N \le X}
     \frac{d(k)}{k}
     \frac{b_N}{\pi^2} 
    \left( \frac{-4N}{k} \right)
     G(N)  = \frac{15}{\zeta(3)} X + O_\varepsilon(X^\varepsilon + XY^{-1/2+\varepsilon}),
    \end{align}
    where the outer sum runs over all positive integers $k \in [1,Y]$ that are perfect squares and where $G(N)$ is defined in \eqref{eq:defnofG}.
\end{lemma}

In our proof of Lemma~\ref{lem:squarecontribution}, we shall 
invoke the P\'{o}lya-Vinogradov inequality in the following form. 
\begin{lemma}[{Pólya--Vinogradov inequality, see \cite[Theorem 9.18]{multntmontvaugh}}] \label{lem:ploya-vino}
    Let $\chi$ be a non-principal character modulo $q$. Then for any integers $M\ge 0, b \ge 0$, we have
    \[
    \sum_{n= 0}^M \chi(n+b) \ll \sqrt{q}\log(q).
    \]
\end{lemma}
Concretely, we use the following direct consequence of Lemma~\ref{lem:ploya-vino}: for all $a \bmod 16$ and all odd integers $q$, with $\overline{16}$ denoting the inverse of $16 \bmod q$, we get
\begin{align}\label{eq:polyavinoapplication16}
    \sum_{\substack{n = 0, \\ n \equiv a \bmod 16}}^M \chi(n) = 
    \sum_{\substack{k=0}}^{\lfloor \frac{M-a}{16} \rfloor} \chi(16k+a)
    = \chi(16) 
    \sum_{\substack{k=0}}^{\lfloor \frac{M-a}{16} \rfloor} \chi(k+\overline{16}a)
    \ll \sqrt{q}\log(q)
\end{align}
for any nonprincipal character $\chi$. 

We also need the following alternative representations of $G(N)$.

\begin{lemma}\label{lem:gpprodusum}
    With $g_N(p)$ defined as in Theorem~\ref{thm:2-4}, we have 
    \[
    G(N) 
    =
      \sum_{\substack{t^2 \mid N, \\ t \text{ odd}}} \frac{2^{\omega(t)}}{t} \prod_{\substack{p \mid t}} \left( 1 - \tfrac{1}{p} \left( \tfrac{-N/t^2}{p} \right) \right)^{-1}
    = \sum_{\substack{t^2 \mid N, \\ t \text{ odd}}} \frac{2^{\omega(t)}}{t} \prod_{\substack{p \mid t, \\ p \nmid N/t^2}} 
    \left( 
    \frac{
     1 + p^{-1}\left( \tfrac{-N/t^2}{p} \right)
    }{1- p^{-2}}
     \right).
    \]
\end{lemma}
\begin{proof}
    The formula for geometric series gives
    \[
        g_N(p) = 2 \sum_{1 \le k \le \nu_p(N)/2} \frac{1}{p^k} \left( 1 - \tfrac{1}{p} \left( \tfrac{-N/p^{2k}}{p} \right) \right)^{-1}.
    \]
    Since $ \left( \tfrac{-N/t^2}{p} \right) = 0 $ if $p \mid N/t^2$, we have 
    \[
    \prod_{\substack{p^2 \mid N, \\ p\neq 2}}(1+ g_{N}(p)) = \sum_{\substack{t^2 \mid N, \\ t \text{ odd}}} \frac{2^{\omega(t)}}{t} \prod_{\substack{p \mid t }} 
    \left( 1 - \tfrac{1}{p} \left( \tfrac{-N/t^2}{p} \right) \right)^{-1} 
    = 
    \sum_{\substack{t^2 \mid N, \\ t \text{ odd}}} \frac{2^{\omega(t)}}{t} \prod_{\substack{p \mid t, \\ p \nmid N/t^2}} 
    \left( 1 - \tfrac{1}{p} \left( \tfrac{-N/t^2}{p} \right) \right)^{-1} 
    .
    \]
    Since $\left( \tfrac{-N/t^2}{p} \right) \in \{\pm 1\}$ for $p \nmid N/t^2$, by the formula for geometric series, we get
    \[
    \left( 1 - \tfrac{1}{p} \left( \tfrac{-N/t^2}{p} \right) \right)^{-1} = \sum_{\substack{j \ge 0, \\ 2 \mid j}}
    p^{-j}
    +
    \sum_{\substack{j \ge 0, \\ 2 \nmid j}}
    p^{-j} 
     \left( \tfrac{-N/t^2}{p} \right)
    = 
    \frac{
     1 + p^{-1}\left( \tfrac{-N/t^2}{p} \right)
    }{1- p^{-2}}. \qedhere
    \]
\end{proof}

\begin{proof}[Proof of Lemma~\ref{lem:squarecontribution}]
  For even $k$, the inner sum vanishes, since then $\left(\frac{4}{k}\right)=0$.
So only the case $k=u^2$, where $u$ is odd, remains. We denote the inner sum by 
\[
S(u^2) =  \sum_{N \le X} 
     \frac{b_N}{\pi^2}
     \left(\frac{-4 N }{ u^2 } \right)
          G(N) = 
           \sum_{\substack{N \le X, \\
            (N,u)=1}} 
     \frac{b_N}{\pi^2}
          G(N).
          \]
 
\medskip

We begin by bounding the error term coming from nonprincipal characters involved in $G(N)$. By Lemma~\ref{lem:gpprodusum}, we have 
 \begin{align*}
G(N)  =
 \sum_{\substack{t^2 \mid N, \\ t \text{ odd}}} \frac{2^{\omega(t)}}{t} \prod_{\substack{p \mid t, \\ p\nmid N/t^2}} \!\!
    \Bigg( 
    \frac{
     1 + p^{-1}\left( \tfrac{-N/t^2}{p} \right)
    }{1- p^{-2}}
     \Bigg) 
     = 
     \sum_{\substack{t^2 \mid N, \\ t \text{ odd}}} \frac{2^{\omega(t)}}{t} \prod_{\substack{p \mid t, \\ p\nmid N/t^2}} 
    (1- p^{-2})^{-1}
     \sum_{\substack{d \mid t }}
     \frac{\mu(d)^2}{d}
     \left( \frac{-N/t^2}{d} \right).
 \end{align*}
  So we rewrite $S(u^2)$  as
 \begin{align*}
    S(u^2)&=     
           \sum_{\substack{N \le X \\  (N,u)=1}}   
      \frac{  b_{N}}{\pi^2}
     \sum_{\substack{t^2 \mid N, \\ t \text{ odd}}} \frac{2^{\omega(t)}}{t} \prod_{\substack{p \mid t, \\ p\nmid N/t^2}} 
    (1- p^{-2})^{-1} 
     \sum_{\substack{d \mid t }}
     \frac{\mu(d)^2}{d}
     \left( \frac{-N/t^2}{d} \right) \\
     &= \sum_{
     \substack{t \le \sqrt{X}, \\ 
      (t,2u)=1
     }
     } 
     \frac{2^{\omega(t)}}{t}\prod_{\substack{p \mid t}} 
     \frac{1}{(1- p^{-2})}
     \sum_{\substack{d \mid t }}
     \frac{\mu(d)^2}{d}
     \sum_{a=0}^{15}
     \frac{b_{a}}{\pi^2} \!\!
     \sum_{ \substack{M \le X/t^2, \\  (M, u)=1 \\
     Mt^2 \equiv a \bmod 16}} \!\!\!\!\!\!
      \left( \frac{-M}{d} \right)
     \prod_{\substack{p \mid t, \\ p\mid M}} 
    (1- p^{-2}).
 \end{align*}
 The innermost sum is 
 \[
    \sum_{ \substack{M \le X/t^2, \\  (M, u)=1 \\
     Mt^2 \equiv a \bmod 16}} 
      \left( \frac{-M}{d} \right)
     \prod_{\substack{p \mid t, \\ p\mid M}} 
    (1- p^{-2}) = 
    \sum_{\substack{r \mid t}} \frac{\mu(r)}{r^2} 
     \left( \frac{-r}{d} \right)
      \sum_{ \substack{M' \le X/(r t^2), \\  (M', u)=1 \\
     M'r t^2 \equiv a \bmod 16}} 
     \left( \frac{M'}{d} \right),
 \]
 and, for squarefree $d\neq 1$, the sum over $M'$ may be rewritten and bounded as 
 \begin{align*}
     \sum_{ \substack{M' \le X/(r t^2), \\  (M', u)=1 \\
     M'r t^2 \equiv a \bmod 16}} 
     \left( \frac{M'}{d} \right) &=  \sum_{ \substack{M' \le X/(r t^2),  \\
     M'r t^2 \equiv a \bmod 16}} 
     \left( \frac{M'}{d} \right) \sum_{\tau \mid  (M', u)} \mu(\tau) \\
     &= \sum_{\tau \mid u} \mu(\tau) 
     \left( \frac{\tau}{d} \right)
      \sum_{ \substack{M'' \le X/(r t^2\tau),  \\
     M''r t^2 \tau \equiv a \bmod 16}} 
     \left( \frac{M''}{d} \right) \ll_\varepsilon u^\varepsilon d^{1/2 +\varepsilon},
 \end{align*}
 where we employed the divisor bound and Lemma~\ref{lem:ploya-vino} (see also \eqref{eq:polyavinoapplication16}).
Applying the divisor bound repeatedly to the contribution of $d \neq 1$ to $S(u^2)$ yields
\begin{align*}
S(u^2) &= 
  \sum_{\substack{t \le \sqrt{X}, \\  (t,2u)=1}} 
  \frac{2^{\omega(t)}}{t}
  \sum_{
  \substack{
  M \le X/t^2, \\
   (M,u)=1
  }
  } 
  \frac{b_{Mt^2}}{\pi^2} \prod_{\substack{p | t, \\ p \nmid M}} (1 - p^{-2})^{-1} + O_\varepsilon \Big(u^\varepsilon \sum_{t \le X} t^{-1+\varepsilon} \Big) \\
  &=
  \sum_{\substack{t \le \sqrt{X}, \\  (t,2u)=1}} 
  \frac{2^{\omega(t)}}{t}
  \sum_{
  \substack{
  M \le X/t^2, \\
   (M,u)=1
  }
  } 
  \frac{b_{Mt^2}}{\pi^2} \prod_{\substack{p | t, \\ p \nmid M}} (1 - p^{-2})^{-1} + O_\varepsilon 
  ( u^\varepsilon X^\varepsilon ).
\end{align*}
This contributes at most
\[
  \ll_\varepsilon X^\varepsilon  \sum_{u \le Y^{1/2}} \frac{d(u^2)}{u^{2-\varepsilon}} \ll X^\varepsilon
\]
for the final sum.

\medskip

Now we deal with the mod 16 condition coming from $b_{t^2 M}$: By the orthogonality relations, we rewrite the above as
 \begin{align*} 
 \frac{1}{16}  &        \sum_{\substack{t \le \sqrt{X}, \\  (t,2u)=1}}     
  \frac{2^{\omega(t)}}{t}    
  \sum_{\substack{M \le X/t^2, \\  (M,u)=1}}
   \prod_{\substack{p \mid t, \\ p \nmid M}}  
 \left( 1-p^{-2}\right)^{-1}
  \sum_{a=0}^{15} 
 \frac{b_a}{\pi^2}
   \sum_{h=0}^{15} e_{16}(h(t^2 M -a))
 \\
 =& 
  \frac{1}{16} 
  \sum_{a=0}^{15} 
 \frac{b_a}{\pi^2}
   \sum_{h=0}^{15} e_{16}(-ha)      
   \sum_{\substack{t \le \sqrt{X}, \\  (t,2u)=1}}     
  \frac{2^{\omega(t)}}{t}      
  \sum_{\substack{M \le X/t^2, \\  (M,u)=1}}
   \prod_{\substack{p \mid t, \\ p \nmid M}}  
 \left( 1-p^{-2}\right)^{-1}
    e_{16}(ht^2 M).
  \end{align*}
  Let us assume that $h \neq 0$ and consider the sum 
  \begin{align*}
    \sum_{\substack{M \le X/t^2 }}
     \mathbf{1}_{[ (M,u)=1]} & e_{16}(ht^2 M)
 \prod_{\substack{p \mid t, \\ p \nmid M}}  
 \left( 1-p^{-2}\right)^{-1} 
 \\
&   =  
\sum_{\substack{M \le X/t^2 }} 
\sum_{q \mid  (M,u)} \mu(q) \ e_{16}(ht^2 M) \prod_{\substack{p \mid t, \\ p \nmid M}}  
 \left( 1-p^{-2}\right)^{-1}  
 \\
&   = 
\prod_{\substack{p \mid t}}  
 \left( 1-p^{-2}\right)^{-1} 
\sum_{q | u} \mu(q) \sum_{\substack{M' \le X/(q t^2) }}     e_{16}(ht^2 M'q)   \prod_{\substack{p| (t, M')}}  
 \left( 1-p^{-2}\right)  
  \end{align*}
 Now, noting that $ (t,M'q) =  (t,M')$ as $ (t,q)=1$, we have
  \begin{align*}
   \sum_{q | u} \mu(q) \!\!\! \sum_{\substack{M' \le X/(q t^2) }}  e_{16}(ht^2 M'q) 
   \prod_{\substack{p \mid  (t,M'q)}}  
 \left( 1-p^{-2}\right) 
 &= \sum_{q | u} \mu(q) \sum_{\substack{M' \le X/(q t^2) }}  e_{16}(ht^2 M'q) 
  \sum_{s \mid  (t, M' q)} \frac{\mu(s)}{s^2} \\
 &= \sum_{q | u} \mu(q)
 \sum_{\substack{s | t
 }} 
 \frac{\mu(s)}{s^2}\sum_{\substack{M'' \le X/(q s t^2) }}  e_{16}(ht^2 M''s q)  \\
 &\ll \sum_{q | u} 
 \sum_{\substack{s | t
 }} 1 \ll_\varepsilon u^\varepsilon t^\varepsilon,
  \end{align*}
where the first upper bound follows from the fact that $t^2s q$ is odd, so $t^2s q h \not \equiv 0 \bmod 16$ and therefore  $\sum_{\substack{M'' \le X/(q s t^2) }}  e_{16}(ht^2 M''s q) \ll 1$. The second upper bound follows from the divisor bound.

As $\prod_{\substack{p \mid t}}  
 \left( 1-p^{-2}\right)^{-1}  \ll 1$, we may bound the contribution from $h \neq 0$ to $S(u^2)$ by
 \begin{align*}
     \frac{1}{16} 
  \sum_{a=0}^{15} 
 \frac{b_a}{\pi^2}
   \sum_{h=1}^{15} e_{16}(-ha) &
   \sum_{\substack{t \le \sqrt{X}, \\  (t,2u)=1 }}     
  \frac{2^{\omega(t)}}{t}    
  \sum_{\substack{M \le X/t^2, \\  (M,u)=1}}
   \prod_{\substack{p \mid t, \\ p \nmid M}}  
 \left( 1-p^{-2}\right)^{-1}
    e_{16}(ht^2 M) \\
 &\ll_\varepsilon \sum_{\substack{t \le \sqrt{X} }}     
  \frac{2^{\omega(t)}}{t} u^\varepsilon t^\varepsilon \ll_\varepsilon u^\varepsilon   \sum_{\substack{t \le \sqrt{X}}} t^{-1+2\varepsilon} \ll_\varepsilon u^\varepsilon X^\varepsilon.
 \end{align*}
 Thus we may bound the whole contribution to $\sum_{u \le \sqrt{Y}} S(u^2) \frac{d(u^2)}{u^2}$ from $h \neq 0$ by 
 \[
    \ll_\varepsilon X^\varepsilon \sum_{u \le \sqrt{Y}} u^{-2+2\varepsilon} \ll X^\varepsilon.
 \]

Now we complete the truncated sums. 
As $\frac{1}{16}\sum_{a=0}^{15} b_a = 70$, we may express the left-hand side of \eqref{eq:squaresum} as 
\begin{align*}
      \frac{70}{\pi^2}     \sum_{
     \substack{u \le \sqrt{Y} \\ u \text{ odd}
     }
     }    
  \frac{ d(u^2) }{u^2} 
          \sum_{\substack{t \le \sqrt{X}, \\  (t, 2u)=1}}  
  \frac{2^{\omega(t)}}{t}    
  \sum_{\substack{M \le (X/t^2), \\  (M, u )=1}} 
 \prod_{\substack{p \mid t, \\ p \nmid M}} 
 \left( 1-p^{-2}\right)^{-1} + O(X^\varepsilon).
\end{align*}

Write $M = sW$ and $t = sy$, where $s =  (M,t)$.  Then we have $ (W,y)=1$ and the double sum becomes 
\begin{align*}
    \sum_{\substack{t \le \sqrt{X}, \\  (t, 2u)=1}}  
  \frac{2^{\omega(t)}}{t} 
  \sum_{\substack{M \le (X/t^2), \\  (M, u )=1}} 
 \prod_{\substack{p \mid t, \\ p \nmid M}} 
 \left( 1-p^{-2}\right)^{-1} 
= 
  \sum_{\substack{ s \le \sqrt{X}, \\
  (s,2u)=1}} 
 \sum_{\substack{
 y \le \sqrt{X}/s, \\  (y,2u) =1
 }}
 \frac{2^{\omega(sy)}}{sy}
   \prod_{\substack{p \mid sy, \\ p \nmid s}}
  \left( 1-p^{-2}\right)^{-1} 
  \sum_{\substack{
 W \le X/(s^3y^2), \\  (W,u y) =1
 }}
  1.
\end{align*}
Using the estimate 
\begin{align}\label{eq:comprimesumasymp}
\sum_{\substack{W \le (X/(s^3y^2)), \\  (W , uy)=1}} 
   1 = \frac{\varphi (uy)}{uy} \frac{X}{s^3y^2} + O_\varepsilon((uy)^\varepsilon);
\end{align}
we note that, as $ (u,y)=1$, we have $\varphi(uy)=\varphi(u)\varphi(y)$. 
We have
$\prod_{\substack{p \mid sy}}
  ({1-p^{-2}})^{-1} \le \zeta(2)$
and
$\prod_{\substack{p \mid s}} 
 ( 1-p^{-2}) \le 1$; thus \eqref{eq:comprimesumasymp} introduces an error term of order 
\begin{align*} 
\sum_{u^2 \le Y} u^{-2+\varepsilon}
    \sum_{\substack{ s \le \sqrt{X}, \\
  (s,2u)=1}}& 
 \sum_{\substack{
 y \le \sqrt{X}/s, \\  (y,2u) =1
 }}
 \frac{2^{\omega(sy)}}{sy}
  (uy)^\varepsilon  \ll_\varepsilon \sum_{u^2 \le Y} u^{-2+2\varepsilon}
    \sum_{\substack{ s \le \sqrt{X}}} 
 s^{-1+\varepsilon}
 \sum_{\substack{
 y \le \sqrt{X}
 }}
  y^{-1+2\varepsilon} 
  \ll_\varepsilon X^{2\varepsilon},
\end{align*}
with the remaining main term given by 
\begin{align*}
    X  \frac{70}{\pi^2}
    \sum_{\substack{u^2 \le Y, \\ u \text{ odd}}}\frac{d(u^2)\varphi(u)}{u^3} 
    \sum_{\substack{ s \le \sqrt{X}, \\
  (s,2u)=1}} 
 \frac{1}{s^4}
 \sum_{\substack{
 y \le \sqrt{X}/s, \\  (y,2u) =1
 }}
 \frac{2^{\omega(sy)}\varphi(y)}{y^4}  \prod_{\substack{p \mid y, \\ p \nmid s}} 
 \left( 1-p^{-2}\right)^{-1}.
\end{align*}
We now extend the sum over $y \le \sqrt{X}/s$ to a sum over all $y \ge 1$. Indeed, we have  $ 2^{\omega(sy)}\varphi(y) y^{-4}
  \prod_{\substack{p \mid y,  p \nmid s}}
  (1-p^{-2})^{-1}  \ll y^{-3+\varepsilon} s^\varepsilon$ and so 
  \begin{align*}
  X 
    \sum_{
     \substack{u \le \sqrt{Y}, \\ u \text{ odd}
     }
     }
    \frac{d(u^2)\varphi(u)}{u^3} 
     \sum_{\substack{ s \le \sqrt{X}, \\
  (s,2u)=1}} 
 s^{-4}
 \sum_{\substack{
 y > \sqrt{X}/s 
 }} 
 &\frac{2^{\omega(sy)}\varphi(y)}{y^4} 
  \prod_{\substack{p \mid y, \\ p \nmid s}}
  \left(1-p^{-2}\right)^{-1} \\
  & \quad \ll_\varepsilon 
  X 
    \sum_{u^2 \le Y} u^{-2+\varepsilon}
    \sum_{\substack{ s \le \sqrt{X} }} 
 s^{-4+\varepsilon}
 \sum_{\substack{
 y > \sqrt{X}/s
 }}
 y^{-3+\varepsilon} \\
 & \quad \ll_\varepsilon  X^{\varepsilon}  
    \sum_{\substack{ s \le \sqrt{X} }} 
 s^{-2+\varepsilon} \ll_\varepsilon X^{\varepsilon} .
  \end{align*}
  We may similarly extend the sum over $s \le \sqrt{X}$ to a sum over all $s \ge 1$ with an error term of 
  $$
  X \sum_{s > \sqrt{X}}s^{-4+\varepsilon} \ll_\varepsilon X^{-1/2+\varepsilon}.
  $$ 
  Finally, we may extend the sum over all odd $ u \le \sqrt{Y}$ to a sum over all odd $u \ge 1$ with an error term of order $$X \sum_{u \ge \sqrt{Y}} u^{-2+\varepsilon} \ll_\varepsilon X Y^{-1/2 + \varepsilon}.$$

 It thus remains to compute the main term, that is, to evaluate 
 \begin{align*}
     X \frac{70}{\pi^2} \sum_{
     \substack{u \ge 1, \\ u \text{ odd}
     }
     } \frac{d(u^2)\varphi(u)}{u^3} \sum_{
     \substack{s \ge 1, \\  (s,2u)=1}
     }
     s^{-4}
 \sum_{\substack{
 y \ge 1, \\  (y,2u)=1 
 }} 
 \frac{2^{\omega(sy)}\varphi(y)}{y^4}
  \prod_{\substack{p \mid y, \\ 
  p\nmid s}}
 \left( 1-p^{-2}\right)^{-1}.
 \end{align*}
 As the function $y \mapsto 2^{\omega(sy) - \omega(s)}\varphi(y) y^{-4}
  \prod_{\substack{p \mid y, 
  p\nmid s}}
 \left( 1-p^{-2}\right)^{-1} $ is multiplicative, we have  
 \begin{align*}
      \sum_{\substack{
 y \ge 1, \\  (y,2u)=1 
 }} \!\!\!\!
 \frac{2^{\omega(sy)}\varphi(y)}{y^4}
  \prod_{\substack{p \mid y, \\ 
  p\nmid s}}
 \left( 1-p^{-2}\right)^{-1} \!
 &= 2^{\omega(s)} \!
    \prod_{\substack{p \nmid 2u, \\ p \mid s}} \!
    \bigg(1+\frac{1}{p(p^2+p+1)}\bigg) \!
    \prod_{\substack{p \nmid 2u ,\\ p \nmid s}} \!
    \bigg(1+\frac{2p}{(p^2-1)(p^2+p+1)}\bigg)
    \\
    &= 
    2^{\omega(s)} \!
\prod_{p\nmid 2u} \!
\left(
1+\frac{2p}{(p^2-1)(p^2+p+1)}
\right) \!
\prod_{\substack{p\nmid 2u\\ p\mid s}}
\frac{(p-1)(p+1)^2}{p(p^2+p-1)}.
 \end{align*}
Similarly, since $s \mapsto s^{-4} 2^{\omega(s)}\prod_{\substack{p\nmid 2u, p\mid s}}
 (p-1)(p+1)^2(p(p^2+p-1))^{-1}$ is multiplicative,
 \begin{align*}
     \sum_{
     \substack{s \ge 1, \\  (s,2u)=1}
     }
     \frac{2^{\omega(s)}}{s^4}
     \prod_{\substack{p\nmid 2u\\ p\mid s}}
\frac{(p-1)(p+1)^2}{p(p^2+p-1)} = \prod_{p\nmid 2u}
\left(
1+
\frac{2(p+1)}
{p(p^2+1)(p^2+p-1)}
\right).
 \end{align*}
 Finally, we have  
 \begin{align*}
     \sum_{
     \substack{u \ge 1, \\ u \text{ odd}
     }
     } \frac{d(u^2)\varphi(u)}{u^3} \prod_{p\nmid 2u}
     &  
\left(
1+
\frac{2(p+1)}
{p(p^2+1)(p^2+p-1)}
\right)    
\left(
1+\frac{2p}{(p^2-1)(p^2+p+1)}
\right) \\
&\qquad\qquad\qquad\qquad\qquad\qquad\qquad\quad= 
 \sum_{
     \substack{u \ge 1, \\ u \text{ odd}
     }
     } \frac{d(u^2)\varphi(u)}{u^3} \prod_{p\nmid 2u}
     \left(
1+\frac{2}{p(p^2-1)}
\right) \\
&\qquad\qquad\qquad\qquad\qquad\qquad\qquad\quad= \prod_{p\neq 2}
\frac{(1-p^{-3})(1-p^{-4})}{(1-p^{-2})^3}
=
\frac{3\pi^2}{14\zeta(3)}.
 \end{align*}
 Together with $\frac{1}{16}\sum_{a=0}^{15}b_a=70$,
this yields the leading constant
\[
\frac{70}{\pi^2}\cdot
\frac{3\pi^2}{14\zeta(3)}
=\frac{15}{\zeta(3)},
\]
as claimed.

\end{proof}
\subsection{Applying Poisson summation to the non-square contribution}\label{section:applying_poisson}
We now apply Poisson summation to the contribution
from non-square $k$ in \eqref{eq:asumdecompwithY}.

\medskip

Throughout the following computations, let $ X \ge Z > 2$ be a smoothing parameter (as specified in \eqref{eq:choiceofyz})
 and let $w: \mathbf{R} \to [0,1]$ be a smooth function that is compactly supported in $(0,1)$ with the properties that $w(x)=1$ for $1/Z \le x \le 1-1/Z$ and $w^{(j)}(x) \ll_j Z^j$ for all $j \ge 1$.\footnote{For concreteness, one may choose $w(x) = F(Zx)F(Z(1-x))$, where we set 
$$F(x) = \frac{\eta(x- \tfrac{1}{2})}{ \eta (x-\tfrac{1}{2}) + \eta (1-x) }, \qquad \text{with} \qquad \eta(x) = e^{-1/x} \mathbf{1}_{(0,\infty)}(x).$$}

We write 
    \begin{equation*}
    % \label{def:defw_hat}
        \widehat{w}(t) = \int_{-\infty}^\infty w(x) e(-tx) \mathrm{d}x
    \end{equation*}
    for the Fourier transform of $w$.
    Since $w$ is smooth on $\R$ and compactly supported, it is bounded on $[0,1]$ and vanishes outside of this interval. Thus it is integrable on $\R$ and 
    \begin{align*}
    |\widehat{w}(t)| \ll \int_{0}^1 |w(x)| \mathrm{d}x \ll 1.
    \end{align*}
   As $w$ is constant on $[1/Z, 1-1/Z]$ and on $\R\setminus [0,1]$, the support of all derivatives $w^{(j)}$ is contained inside $[0,1/Z] \cup [1-1/Z,1]$. Repeated integration by parts yields
   \begin{align*}
       |\widehat{w}(t)| \ll_j \lambda\left( [0,\tfrac{1}{Z}]
       \cup  [1-\tfrac{1}{Z},1]\right)\frac{Z^j}{|t|^j} \ll_j \frac{Z^{j-1}}{|t|^j}
   \end{align*}
   for all $j\ge 1$, where $\lambda$ denotes the Lebesgue measure.
   Therefore we have 
   \begin{align}\label{eq:fouriertrafobound}
        |\widehat{w}(t)| \ll_j \min \left\{ 1, \frac{Z^{j-1}}{|t|^j}\right\}.
   \end{align}
   Similarly, we may bound the $k$-th derivative as  
   \begin{align}\label{eq:whatderivativebound}
       \left|\widehat{w}^{(k)}(t)\right| \ll_{k,j} \min \left\{ 1, \frac{Z^{j-1}}{|t|^j}\right\}.
   \end{align}
The main goal of this section is to prove the following lemma.

\begin{lemma}\label{lem:possionsummationapplication}
    We have 
    \begin{align*}
      \sum_{ 
     \substack{k \le Y, \\ k \neq \square
     }
     } 
     \sum_{N \le X}
     \frac{d(k)}{k}
     \frac{b_N}{\pi^2} 
    \left( \frac{-4N}{k} \right)
     G(N) =& X \sum_{
     \substack{k \le Y, \\ k \neq \square, \\ k \text{ odd}
     }
     }
     \frac{d(k)}{k} \!\!\!
     \sum_{\substack{t \le \sqrt{X}, \\  (t,2k)=1}} \! \! \! \frac{2^{\omega(t)}}{t^3} 
         \sum_{r \in \Z}
         \widehat{w}\big(\tfrac{X r}{Q_{k,t}t^2} \big)
         H_{k,t}(r) \\
         &+O_\varepsilon(
         (XYZ)^\varepsilon X^{1/2} Z^{-1/4}
         ).
    \end{align*}
    Here we set 
      $Q_{k,t}=16\rad(kt)$ and 
   \[
    H_{k,t}(r) = \frac{1}{Q_{k,t}} \sum_{q \bmod Q_{k,t}}\frac{b_{qt^2}}{\pi^2} 
    \left( \frac{-q}{k} \right)
    \prod_{\substack{p \mid t}} \left( 1 - \tfrac{1}{p} \left( \tfrac{-q}{p} \right) \right)^{-1} e_{Q_{k,t}}(qr),
   \]
   where $b_N$ is defined in Theorem~\ref{thm:2-4}.
\end{lemma}

We will compute $H_{k,t}(r)$ further as follows.

\begin{lemma}\label{lem:evaluationqsum}
    Let $k \ge 1$ be an odd integer that is not a perfect square, let $t$ be a positive integer with $ (2k,t)=1$, let $r \in \Z$ and let $H_{k,t}$ and $Q_{k,t}$ be defined as in Lemma~\ref{lem:possionsummationapplication}.
    Further, we factor $k=k_1 k_2^2$, 
    where 
    \[
    k_1 = \prod_{2 \nmid \nu_p(k)} p^{\nu_p(k)}, \qquad k_2 = \prod_{2 \mid \nu_p(k)} p^{\nu_p(k)/2}.
    \]
    Write $\kappa = \rad(k), \kappa_1 = \rad(k_1), \kappa_2=\rad(k_2)$ and $\tau = \rad(t)$.
    We have
    \begin{align*}
    H_{k,t}(r) = A_{k,t}(r) B_{k,t}(r) C_{k_2}(r) D_{k,t}(r),
    \end{align*}
    where 
    \begin{align*}
         A_{k,t}(r) &= \frac{1}{16}\sum_{a \bmod 16}  \frac{b_{a}}{\pi^2} e_{16}(ra \kappa \tau)
         \\ B_{k,t}(r) &= \left(\frac{-r \kappa_2 \tau}{\kappa_1 }
    \right)
     \frac{\epsilon_{\kappa_1}}{\sqrt{\kappa_1}}, \\
     C_{k_2}(r) &=  \frac{\mu(\kappa_2) }{\mu( (\kappa_2,r))}
    \frac{ \varphi( (\kappa_2,r))}{\kappa_2}, \\ D_{k,t}(r)&= 
          \prod_{p \mid t} (p(1-p^2))^{-1} \prod_{\substack{p \mid t, \\ p \mid r}} (1-p^3)
          \sum_{\substack{x \mid \tau 
          }}
          \mu(x) x^{3/2}
          \left( \frac{-r \kappa \tau/x}{x}
          \right)\epsilon_x,
    \end{align*}
   where, for any odd integer $\ell$
   \[
    \epsilon_\ell = \begin{cases}
        1 & \text{ if $\ell \equiv 1 \bmod 4$}, \\
        i  & \text{ if $\ell \equiv 3 \bmod 4$}.
    \end{cases}
    \]
\end{lemma}

For the proof of Lemma~\ref{lem:possionsummationapplication}, we need the following large sieve inequality due to Heath-Brown. 
 \begin{lemma}[{\cite[Corollary~2]{hbinequ}}] \label{lem:hbrealchar}
    Let $N, Q$ be positive integers and let $a_1, \ldots, a_N$ be arbitrary complex numbers. Then we have
\[
\sum_{d\le Q} \mu(2d)^2  \bigg| \sum_{n \leq N}  a_n \left( \frac{n}{d} \right) \bigg|^2 
\;\;\ll_\varepsilon\;\; (QN)^\varepsilon  (Q + N) \sum_{n_1n_2 = \square} |a_{n_1}a_{n_2}|,
\]
for any $\varepsilon > 0$.
\end{lemma}

\begin{proof}[Proof of Lemma~\ref{lem:possionsummationapplication}]

    For any odd integer $k \ge 1$, we define the smooth sum
    \[
    F_{w}(k) =  
     \sum_{N \ge 1}
     \frac{b_N}{\pi^2} 
    \left( \frac{-N}{k} \right)
     G(N) w(X^{-1}N).
    \]
     Since $w$ is compactly supported in $(0,1)$, only finitely many summands are non-zero.
    By Lemma~\ref{lem:gpprodusum}, we have 
    \begin{align*}
     F_{w}(k) &= \sum_{N \ge 1}
        \frac{b_{N}}{\pi^2} 
    \left( \frac{-N}{k} \right)
    \sum_{\substack{t^2 \mid N, \\ t \text{ odd}}} \frac{2^{\omega(t)}}{t} \prod_{\substack{p \mid t}} \left( 1 - \tfrac{1}{p} \left( \tfrac{-N/t^2}{p} \right) \right)^{-1}
      w(X^{-1}N) \\
      &=\sum_{\substack{t \le \sqrt{X}, \\  (t,2k)=1}} \frac{2^{\omega(t)}}{t} \sum_{M \in \Z} \frac{b_{Mt^2}}{\pi^2} 
    \left( \frac{-M}{k} \right)
    \prod_{\substack{p \mid t}} \left( 1 - \tfrac{1}{p} \left( \tfrac{-M}{p} \right) \right)^{-1}
      w(X^{-1}Mt^2).
    \end{align*}
    We observe that for fixed $t$, the function of $M$ given by 
    \begin{equation}\label{eq:deflowercasehkt}
          h_{k,t}(M) = \frac{b_{Mt^2}}{\pi^2} 
    \left( \frac{-M}{k} \right)
    \prod_{\substack{p \mid t}} \left( 1 - \tfrac{1}{p} \left( \tfrac{-M}{p} \right) \right)^{-1}  
    \end{equation}
    is periodic with a period of $Q_{k,t} = 16\prod_{p\mid k t}p$. Therefore, by the orthogonality of characters, we have 
    \begin{align*}
         F_{w}(k) 
         &= \sum_{\substack{t \le  \sqrt{X}, \\  (t,2k)=1}}  \frac{2^{\omega(t)}}{tQ_{k,t}} 
         \sum_{M \in \Z} w(X^{-1}Mt^2)
         \sum_{q_1,q_2 \bmod Q_{k,t}}
         h_{k,t}(q_1)e_{Q_{k,t}}(q_2(M-q_1)) \\
         &=
         \sum_{\substack{t \le \sqrt{X}, \\  (t,2k)=1}}   \frac{2^{\omega(t)}}{tQ_{k,t}}  
         \sum_{q_1,q_2 \bmod Q_{k,t}} 
         h_{k,t}(q_1)
         e_{Q_{k,t}}(-q_1q_2) 
         \sum_{M \in \Z}  w(X^{-1} Mt^2) e_{Q_{k,t}}(q_2M).
    \end{align*}
   Applying Poisson summation to the innermost sum yields
   \[
   \sum_{M \in \Z} w(X^{-1}Mt^2) e_{Q_{k,t}}(q_2M) = 
   \frac{X}{t^2}
   \sum_{m \in \Z} \widehat{w}(X(Q_{k,t}t^2)^{-1}(Q_{k,t}m-q_2)).
   \]
Since $r = Q_{k,t}m-q_2$ takes every integer value exactly once as $m$ runs through the integers and $q_2$ runs
 through the different residue classes $\bmod~Q_{k,t}$, we have 
   \[
   \sum_{q_2 \bmod Q_{k,t}} 
   \sum_{m \in \Z} \widehat{w}(X(Q_{k,t}t^2)^{-1}(Q_{k,t}m-q_2)) e_{Q_{k,t}}(-q_1q_2)
   = 
   \sum_{r \in \Z} \widehat{w}(X(Q_{k,t}t^2)^{-1} r)e_{Q_{k,t}}(q_1r),
   \]
   and therefore 
   \begin{align*}
         F_{w}(k) 
         &= X \sum_{\substack{t \le \sqrt{X}, \\  (t,2k)=1}}  \frac{2^{\omega(t)}}{t^3} 
         \sum_{r \in \Z}
         \widehat{w}(X(Q_{k,t}t^2)^{-1} r)
         H_{k,t}(r) .
   \end{align*}   

   It remains to give an upper bound for the error introduced by smoothing, i.e. for
   \[
   E_{w}(X,Y,Z) := 
    \sum_{ 
     \substack{k \le Y, \\ k \neq \square
     }
     }  
     \sum_{N \le X}
     \frac{d(k)}{k}
     \frac{b_N}{\pi^2} 
    \left( \frac{-4N}{k} \right)
     G(N) \left(
    1 - w (N/X)
     \right).
   \]
   To this end, it suffices to bound 
     \begin{equation}\label{eq:toboundsmoothingerror}
           \sum_{\substack{k \le Y, \\ k \neq \square}}
           \bigg|
     \sum_{N \in I}
     \frac{d(k)}{k}
     \frac{b_N}{\pi^2} 
    \left( \frac{-4N}{k} \right)
     G(N) \bigg| ,
     \end{equation}
     where $I \subset [1,X]$ is a fixed interval of length $|I| \ll XZ^{-1}$.
   Indeed, $1-w(N/X) = 0$ for all $N \in [XZ^{-1},X - X Z^{-1}]$ and therefore $|E_w(X,Y,Z)|$ is bounded from above by 
   \[
     \Bigg| \sum_{ 
     \substack{k \le Y, \\ k \neq \square
     }
     }  
     \sum_{N \in \mathcal{I}_{X,Z}} \!\!
     \frac{d(k)}{k}
     \frac{b_N}{\pi^2} 
    \left( \frac{-4N}{k} \right)
     G(N)
     \Bigg| \! \,
     + \! \,
     \Bigg| \sum_{ 
     \substack{k \le Y, \\ k \neq \square
     }
     }  
     \sum_{N \in \mathcal{I}_{X,Z}} \!\!
     \frac{d(k)}{k}
     \frac{b_N}{\pi^2} 
    \left( \frac{-4N}{k} \right)
     G(N) w(N/X)
     \Bigg|,
   \]
   where $ \mathcal{I}_{X,Z} = [1, XZ^{-1}] \cup 
   [X - X Z^{-1}, X]$.
   Evaluating the inner sum over the interval $[1, XZ^{-1}]$ by partial summation, we have 
   \begin{align*}
          \sum_{N \le XZ^{-1}} 
     \frac{b_N}{\pi^2} 
    \left( \frac{-4N}{k} \right)
       G(N) 
     w (N/X) 
     &\ll 
          \sum_{N \le XZ^{-1}}
     \frac{b_N}{\pi^2} 
    \left( \frac{-4N}{k} \right)
     G(N) \\
     &+
   X^{-1}
   \int_1^{XZ^{-1}}
   |w'(t/X)|\bigg|
          \sum_{N \le t}
     \frac{b_N}{\pi^2} 
    \left( \frac{-4N}{k} \right)
     G(N) \bigg|
     \mathrm{d}t.
   \end{align*}
    Since $|w'(t/X)| \ll Z$ and thus $\int_{1}^{XZ^{-1}} X^{-1} |w'(t/X)| \mathrm{d}t \ll XZ^{-1} X^{-1} Z \ll 1$, we get 
    \begin{align*}
    \sum_{\substack{k \le Y, \\ k \neq \square}} \frac{d(k)}{k} 
     X^{-1}\int_1^{XZ^{-1}}
   w'(t/X) &
     \bigg|
          \sum_{N \le t}
     \frac{b_N}{\pi^2} 
    \left( \frac{-4N}{k} \right)
     G(N) \bigg|
     \mathrm{d}t
     \\
   & \ll 
     \max_{1 \le t 
     \le XZ^{-1}}
       \sum_{\substack{k \le Y, \\ k \neq \square}} \frac{d(k)}{k}
     \bigg| \sum_{ N \le t}
     \frac{b_N}{\pi^2} 
    \left( \frac{-4N}{k} \right)
     G(N) 
     \bigg| 
    \end{align*} 
    and therefore 
    \begin{align*}
     \sum_{\substack{k \le Y, \\ k \neq \square}} \frac{d(k)}{k} \!\!\!\!
         \sum_{N \le XZ^{-1}} 
     \frac{b_N}{\pi^2} 
    \left( \frac{-4N}{k} \right) 
       G(N) 
     w (N/X)
 \ll 
    \max_{\substack{I \subset [1,X], \\
     |I| \le XZ^{-1}}}
       \sum_{\substack{k \le Y, \\ k \neq \square}} \frac{d(k)}{k}
     \bigg| \sum_{N \in I}
     \frac{b_N}{\pi^2} 
    \left( \frac{-4N}{k} \right)
     G(N) 
     \bigg|.
    \end{align*}
    Analogously, for the sum over $[X - X Z^{-1}, X]$ we have 
    \begin{align*}
     \sum_{\substack{k \le Y, \\ k \neq \square}} \frac{d(k)}{k} \!\!\!\!
        \sum_{X- X Z^{-1} \le N \le X} \!\!
     \frac{b_N}{\pi^2} 
    \left( \frac{-4N}{k} \right) 
       G(N) 
     w (N/X)  \ll \max_{\substack{I \subset [1,X], \\
     |I| \le X Z^{-1} }}
       \sum_{\substack{k \le Y, \\ k \neq \square}} \frac{d(k)}{k}
     \bigg| \sum_{N \in I}
     \frac{b_N}{\pi^2} 
    \left( \frac{-4N}{k} \right)
     G(N) 
     \bigg|.
    \end{align*}
    
    We estimate \eqref{eq:toboundsmoothingerror}
by treating small and large values of $k$
separately, using Pólya--Vinogradov (in the form of \eqref{eq:polyavinoapplication16}) for the
former and Heath-Brown's quadratic large
sieve (Lemma~\ref{lem:hbrealchar}) for the latter.

    We begin with an upper bound for the inner sum in \eqref{eq:toboundsmoothingerror} that is useful for small $k$. By Lemma~\ref{lem:gpprodusum}, we have 
    \begin{align*}
    \sum_{N \in I}
     \frac{b_N}{\pi^2} 
    \left( \frac{-4N}{k} \right)
           G(N) 
           &= \!  \sum_{N \in I}
     \frac{b_N}{\pi^2} 
    \left( \frac{-4N}{k} \right)\sum_{\substack{t^2 \mid N, \\ t \text{ odd}}} \frac{2^{\omega(t)}}{t} \! \prod_{\substack{p \mid t, \\ p \nmid N/t^2}}
    \left(
    \frac{
     1 + p^{-1}\left( \tfrac{-N/t^2}{p} \right)
    }{1- p^{-2}}
    \right)\\
    &= \!
    \sum_{\substack{t \le \sqrt{X}, \\ t \text{ odd}}} \! \frac{2^{\omega(t)}}{t}  \!
    \left(
    \frac{4t^2}{k}
    \right) \!
    \sum_{t^2 M \in I} \frac{b_{t^2M}}{\pi^2}
    \left(
    \frac{-M}{k}
    \right)
    \prod_{\substack{p \mid t, \\ p \nmid M}} 
    \left(
    \frac{
     1 + p^{-1}\left( \tfrac{-M}{p} \right)
    }{1- p^{-2}}
    \right).
    \end{align*}
    Now we have 
    \begin{align*}
        \sum_{t^2 M \in I} \frac{b_{t^2M}}{\pi^2}
    &\left(
    \frac{-M}{k}
    \right) 
    \prod_{\substack{p \mid t, \\ p \nmid M}} 
    \left(
    \frac{
     1 + p^{-1}\left( \tfrac{-M}{p} \right)
    }{1- p^{-2}}
    \right)  \\
    &= 
    \prod_{\substack{p \mid t}} 
    \left(
    1- p^{-2} \right)^{-1}
    \sum_{t^2 M \in I} \frac{b_{t^2M}}{\pi^2}
    \left(
    \frac{-M}{k}
    \right) 
    \prod_{\substack{p \mid t,  \\ p \mid M}}
    \left(
    1- p^{-2} \right)
    \prod_{\substack{p \mid t}} 
    \left(
     1 + p^{-1}\left( \tfrac{-M}{p} \right) 
    \right)
    \end{align*} 
    and
 \begin{align*}
    \sum_{t^2 M \in I} \frac{b_{t^2M}}{\pi^2} 
    \sum_{r \mid  (t,M)} \frac{\mu(r)}{r^2} \sum_{s \mid t} \frac{\mu(s)^2}{s} \left( 
    \frac{-M}{ks}
    \right) 
    = 
     \sum_{r \mid t}
     \frac{\mu(r)}{r^2}
     \sum_{s \mid t} 
      \frac{\mu(s)^2}{s} 
       \left( 
    \frac{r}{ks}
    \right)
    \sum_{t^2 r M' \in I} \frac{b_{t^2M'r}}{\pi^2} 
    \left( 
    \frac{-M'}{ks}
    \right).
    \end{align*} 
    Let us note that $ (t,k)=1$, so $ (k,s)=1$; this means that, since $k$ is not a square, neither is $ks$. 
    By Lemma~\ref{lem:ploya-vino}, splitting modulo 16 to account for $b_{t^2M'r}$ and employing \eqref{eq:polyavinoapplication16}, we have 
\[
    \sum_{t^2 r M' \in I} \frac{b_{t^2M'r}}{\pi^2} 
    \left( 
    \frac{-M'}{ks}
    \right) \ll_\varepsilon (ks)^{1/2+\varepsilon}.
\]
This gives by the divisor bound that
\[
\sum_{t^2 M \in I} \frac{b_{t^2M}}{\pi^2}
    \left(
    \frac{-M}{k}
    \right) 
    \prod_{\substack{p \mid t, \\ p \nmid M}} 
    \left(
    \frac{
     1 + p^{-1}\left( \tfrac{-M}{p} \right)
    }{1- p^{-2}}
    \right) \ll_\varepsilon t^\varepsilon k^{1/2+\varepsilon}
\] and hence 
\begin{equation*}
     \sum_{N \in I}
     \frac{b_N}{\pi^2} 
    \left( \frac{-4N}{k} \right)
           G(N)  \ll_\varepsilon \sum_{t \le \sqrt{X}} t^{-1+\varepsilon} k^{1/2+\varepsilon}
           \ll k^{1/2+\varepsilon} X^\varepsilon.
\end{equation*}
For $\widetilde{Y} \ge 1$, we have 
\begin{align}\label{eq:smoothingerror_smallk}
    \sum_{\substack{ k \le \widetilde{Y} , \\ k \neq \square}}
     \frac{d(k)}{k}
     \left|
     \sum_{N \in I}
     \frac{b_N}{\pi^2} 
    \left( \frac{-4N}{k} \right)
           G(N)\right|  \ll_\varepsilon
           X^\varepsilon 
           \sum_{\substack{ k \le \widetilde{Y} , \\ k \neq \square}}
           k^{-1/2 + \varepsilon} \ll_\varepsilon
            X^\varepsilon \widetilde{Y}^{1/2+\varepsilon}.  
\end{align}

 We next bound the contribution from large $k$. We have 
 \begin{align*}
 \sum_{\substack{K \le k \le 2K, \\ k \neq \square}}
     \frac{d(k)}{k}
     \left|
     \sum_{N \in I}
     \frac{b_N}{\pi^2} 
    \left( \frac{-4N}{k} \right)
           G(N) \right|
           &\ll_\varepsilon K^{-1+\varepsilon}
            \sum_{\substack{K \le k \le 2K, \\ k \neq \square }}
     \Bigg|
     \sum_{N \in I}
     \frac{b_N}{\pi^2} 
    \left( \frac{4N}{k} \right)
           G(N)
           \Bigg|.
 \end{align*}
 Noting that $\left( \frac{-4N}{k} \right) = 0$ for even $k$, we may restrict to odd numbers $k$. We may thus decompose $k = rm^2$ for some odd squarefree number $r>1$; this gives 
 \begin{align*}
     \sum_{\substack{K \le k \le 2K, \\ k \neq \square }}
     \Bigg|
     \sum_{N \in I}
     \frac{b_N}{\pi^2} 
    \left( \frac{N}{k} \right)
           G(N)
           \Bigg| 
           &\ll \!\! \! \sum_{m \le (2K)^{1/2}} \sum_{1<r \le 2K/m^2} \mu(2r)^2 \bigg|
           \sum_{N \in I}
     \frac{b_N}{\pi^2} 
      \left( \frac{N}{rm^2} \right)
           G(N)
           \bigg| \\
           &\ll \!\! \!
           \sum_{m \le (2K)^{1/2}}\!\!
           \frac{K^{1/2}}{m}
           \left(
           \sum_{1<r \le 2K/m^2} \mu(2r)^2 \bigg|
           \sum_{N \in I}
     \frac{b_N}{\pi^2} 
      \left( \frac{N}{rm^2} \right)
           G(N)
           \bigg|^2
           \right)^{1/2} \!,
 \end{align*}
 where we employed  Cauchy's inequality. For the inner sum, we use  Lemma~\ref{lem:hbrealchar} with $a_N = \frac{b_N}{\pi^2} \left( \tfrac{N}{m^2}
 \right) G(N) \mathbf{1}_I (N)$, noting that $G(N)\ll_\varepsilon N^\varepsilon \ll X^{\varepsilon}$. This yields 
 \begin{align*}
     \sum_{r \le 2K/m^2} \! \mu(2r)^2 \bigg|
           \sum_{N \le X}
     \! a_N \!
      \left( \frac{N}{r} \right)
           \bigg|^2 
        &\ll_\varepsilon (KX)^\varepsilon \left(\frac{K}{m^2} + X\right) X^{\varepsilon}\# \{ (n_1, n_2) \in I^2: n_1 n_2 =\square \}.
 \end{align*}
 
 The number of tuples of integers $(n_1, n_2) \in I^2$ such that $n_1 n_2$ is a square is $\ll X^{1+\varepsilon} Z^{-1}$. This may be proved along the same lines as in the proof of \cite[Corollary~3]{hbinequ}; we include the argument for completeness.
 
 Writing $n_1 = m_1^2 d_1$, $n_2 = d_2m_2^2$ with $d_1, d_2$ squarefree, the product $n_1 n_2$ is a square if and only if $d_1 = d_2 =:d$.  There are $\ll XZ^{-1}$ pairs $(n_1,n_2)$
with $n_1=n_2$. Assume that $n_1 \neq n_2$ and set $I = [I_1,I_2]$. For any integer $d$, there exist distinct integers $m_1 \neq m_2$ with $dm_1^2, dm_2^2 \in I$ only if 
\[
    \frac{I_2^{1/2}}{d^{1/2}} -
    \frac{I_1^{1/2}}{d^{1/2}} \ge 1,
\]
which implies $d \le (I_2^{1/2} - I_1^{1/2})^2 \le (I_2^{1/2} - I_1^{1/2})(I_2^{1/2} + I_1^{1/2}) = I_2 - I_1 \le XZ^{-1}$. This yields 
\begin{align*}
    \# \{ (n_1, n_2) \in I^2: n_1 n_2 =\square \} 
    \ll \# I + \sum_{d \le XZ^{-1}} \left(
     \frac{I_2^{1/2}}{d^{1/2}} -
    \frac{I_1^{1/2}}{d^{1/2}}
    \right)^2 
    &\ll 
    XZ^{-1} + XZ^{-1} \sum_{d \le XZ^{-1}} \frac{1}{d} \\
    &\ll X^{1+\varepsilon}Z^{-1} . 
\end{align*}
It follows that 
\begin{align*}
     \sum_{r \le 2K/m^2} \mu(2r)^2 \bigg|
           \sum_{N \le X}
    a_N
      \left( \frac{N}{r} \right)
           \bigg|^2 \ll_\varepsilon (KX)^\varepsilon \left(\frac{K}{m^2} + X\right)  X^{1+\varepsilon}Z^{-1},
\end{align*}
and thus 
\begin{align*}
     \sum_{\substack{K \le k \le 2K, \\ k \neq \square}}
     \frac{d(k)}{k}\left|
     \sum_{N \in I}
     \frac{b_N}{\pi^2} 
    \left( \frac{-4N}{k} \right)
           G(N) \right|
    &\ll_\varepsilon 
    (KX)^\varepsilon \sum_{m \le (2K)^{1/2}} 
     \left(\frac{X^{1/2}}{Z^{1/2}m^2} + \frac{X}{m(ZK)^{1/2}}\right)  \\
     &\ll_\varepsilon 
    (KX)^{2\varepsilon} 
     \left(\frac{X^{1/2}}{Z^{1/2}} + \frac{X}{(ZK)^{1/2}}\right)  
\end{align*}
Summing dyadically finally gives the bound
\begin{align}\label{eq:largeksmootingerror}
    \sum_{\substack{\widetilde{Y} \le k \le Y, \\ k \neq \square}}
     \frac{d(k)}{k}
     \left|
     \sum_{N \in I}
     \frac{b_N}{\pi^2} 
    \left( \frac{-4N}{k} \right)
           G(N)\right|  \ll_\varepsilon
           (YX)^{\varepsilon} 
     \left(\frac{X^{1/2}}{Z^{1/2}} + \frac{X}{(Z\widetilde{Y})^{1/2}}\right) 
\end{align}

If $Y<X/Z^{1/2}$, the estimate
\eqref{eq:smoothingerror_smallk} alone gives
the required bound with $\widetilde{Y} = Y$. Otherwise, we choose
$\widetilde{Y}=X/Z^{1/2}$.
By \eqref{eq:smoothingerror_smallk} and
\eqref{eq:largeksmootingerror}, this gives 
\[
    \sum_{\substack{ k \le Y, \\ k \neq \square}}
    \frac{d(k)}{k}
    \left|
     \sum_{N \in I}
     \frac{b_N}{\pi^2} 
    \left( \frac{-4N}{k} \right)
           G(N) \right|
           \ll_\varepsilon (XYZ)^\varepsilon 
           X^{1/2} Z^{-1/4}
           . \qedhere
\]
\end{proof}

\begin{proof}[Proof of Lemma~\ref{lem:evaluationqsum}]
    We have, with $h_{k,t}$ as in \eqref{eq:deflowercasehkt}, that 
    \[
     h_{k,t}(q) = \frac{b_{qt^2}}{\pi^2} \left(\frac{-q}{k}\right)
     \prod_{p \mid t} 
     \left(
     1 - p^{-1}
     \left(
     \tfrac{-q}{p}
     \right)
     \right)^{-1},
     \]
     where $b_{qt^2}$ only depends on $q \bmod 16$, $\left(\frac{-q}{k}\right)$ only depends on $q \bmod \kappa_1 \kappa_2$ and each factor $1/(
     1 - p^{-1}
     (
     \tfrac{-q}{p}
     )
     )$ only depends on $ q \bmod p$.
 
     Since the numbers 16, $\kappa_1$, $\kappa_2$ and $\tau$ are pairwise coprime, for $m \in \{p: p \mid k_2  \tau\} \cup \{ 16,  \kappa_1 \}$ there exist by the Chinese remainder theorem integers $u_{m}$ such that $u_m R_{k,t}(m) \equiv 1\bmod m$ with $R_{k,t}(m) = Q_{k,t}/m$ and  
   
    \[
     \frac{1}{Q_{k,t}}
         \sum_{q \bmod Q_{k,t}} h_{k,t}(q)e_{Q_{k,t}}(qr) =
         A_{k,t}(r)
         B_{k,t}(r)
         \prod_{p \mid k_2}
         C_{p}(r) 
         \prod_{p \mid t} D_{p}(r),
    \]
    where 
    \begin{align*}
    A_{k,t}(r) = \frac{1}{16}\sum_{a \bmod 16}  \frac{b_{t^2a}}{\pi^2} e_{16}(ru_{16}a), \quad B_{k,t}(r) = \frac{1}{\kappa_1}\sum_{b \bmod \kappa_1} 
    \left(
    \frac{-b}{\kappa_1}
    \right) e_{\kappa_1}(ru_{\kappa_1}b)
    \end{align*}
    and 
    \begin{align*}
    C_{p}(r) = \frac{1}{p}\sum_{c \bmod p} 
    \left(
    \frac{-c}{p}
    \right)^2 e_{p}(ru_{p}c), \quad 
    D_{p}(r)=
        \frac{1}{p} \sum_{d \bmod p}
        \left(
        1 - p^{-1} 
        \left(
        \tfrac{-d}{p}
        \right)
        \right)^{-1}  e_{p}(ru_{p}d).
    \end{align*}

    First we evaluate $A_{k,t}(r)$.
    As $t$ is odd, we have $t^2 \in \{1, 9\} \bmod 16.$ By the definition of $b_N$ in Theorem~\ref{thm:2-4}, we get $b_a = b_{9a}$ for all $a \bmod 16$ and so 
    \[
    A_{k,t}(r) = \frac{1}{16}\sum_{a \bmod 16}  \frac{b_{t^2a}}{\pi^2} e_{16}(ru_{16}a) = \frac{1}{16}\sum_{a \bmod 16}  \frac{b_{a}}{\pi^2} e_{16}(ru_{16}a).
    \]
    Furthermore, as $\kappa  \tau$ is odd, multiplication by $(\kappa \tau)^2$ is a bijection modulo 16. Thus we may apply the change of variables $a = (\kappa \tau)^2 a'$; again we get $b_{a} = b_{a'(\kappa \tau)^2} = b_{a'}$. Furthermore, we have $u_{16}( \kappa \tau)^2 \equiv \kappa \tau \bmod 16$ and hence
    \[
    A_{k,t}(r) = \frac{1}{16}\sum_{a' \bmod 16}  \frac{b_{a'}}{\pi^2} e_{16}(ra'\kappa \tau).
    \]

Now we evaluate $B_{k, t}(r)$. Since $\kappa_1$ is odd and squarefree,
the quadratic character modulo $\kappa_1$
is primitive, and so the classical evaluation of the Gauss sum yields
\[
    B_{k, t}(r) = \frac{\epsilon_{\kappa_1}}{\sqrt{\kappa_1}} 
    \left(
    \frac{-ru_{\kappa_1}}{\kappa_1}
    \right).
\]
As $u_{\kappa_1}^{-1} \equiv 16 \kappa_2 \tau \bmod \kappa_1$, we further observe 
\begin{align*}
    B_{k, t}(r) = \frac{\epsilon_{\kappa_1}}{\sqrt{\kappa_1}} 
    \left(
    \frac{-ru_{\kappa_1}}{\kappa_1}
    \right) = \frac{\epsilon_{\kappa_1}}{\sqrt{\kappa_1}} 
    \left(
    \frac{-ru_{\kappa_1}^{-1}}{\kappa_1}
    \right)
    = 
    \frac{\epsilon_{\kappa_1}}{\sqrt{\kappa_1}} 
    \left(
    \frac{-r}{\kappa_1}
    \right)
    \left(
    \frac{\kappa_2 \tau}{\kappa_1}
    \right)
    .
\end{align*}

Now we compute $\prod_{p \mid k_2} C_{p}(r)$.  By orthogonality, we have 
$
    C_{p}(r) = \mathbf{1}_{p\mid r}-\frac{1}{p}.
$
Hence we get 
\begin{align*}
    \prod_{p \mid k_2} C_{p}(r) =  \prod_{\substack{p \mid k_2, \\ p \nmid r}} (-p^{-1}) \prod_{\substack{p \mid k_2, \\ p \mid r}} 
    (1 - p^{-1})
    =   \frac{\mu(\kappa_2) }{\kappa_2 } 
    \frac{ (\kappa_2,r)}{\mu( (\kappa_2,r))} 
    \prod_{\substack{p \mid k_2, \\ p \mid r}} 
    (1 - p^{-1}) 
    =
    \frac{\mu(\kappa_2) }{\mu( (\kappa_2,r))}
    \frac{ \varphi( (\kappa_2,r))}{\kappa_2}.
\end{align*}

Finally, we evaluate $\prod_{p \mid t}D_p(r)$; we have 
\begin{align*}
      D_{p}(r)
      &=
        \frac{1}{p} \sum_{d \bmod p}
        \left(
        1 - p^{-1} 
        \left(
        \tfrac{-d}{p}
        \right)
        \right)^{-1}  e_{p}(ru_{p}d) \\
        &=
        \frac{1}{p} \sum_{d \bmod p} \bigg(
        1 + \sum_{j\ge 1} p^{-j} \left(
        \frac{-d}{p}
        \right)^j
        \bigg) e_{p}(ru_{p}d) \\
        &=
         \frac{1}{p} \sum_{d \bmod p} \bigg(
        1 +  \frac{1}{p(1-p^{-2})}  \left(
        \frac{-d}{p}
        \right)+ 
        \frac{ 1}{p^2(1-p^{-2})} \left(
        \frac{-d}{p}
        \right)^2
        \bigg) e_{p}(ru_{p}d) \\ 
        &= 
        \left(
        1+ \frac{1}{p^2(1-p^{-2})}
        \right) \mathbf{1}_{[p \mid r]} -  \frac{1}{p^3(1-p^{-2})} + \frac{\epsilon_p}{p^{3/2}(1-p^{-2})} 
        \left(
        \frac{-ru_p}{p}
        \right);
\end{align*}
this gives 
\begin{align*}
     D_{p}(r) 
     = \begin{cases}
        (p(1-p^2))^{-1}
         (1-p^3)& \text{ if $p \mid r$}, \\[8pt]
         (p(1-p^2))^{-1} \left(1-
         \epsilon_p \left( 
         \frac{-ru_p}{p}
         \right)p^{3/2} 
         \right)
         & \text{ if $p \nmid r$}.
     \end{cases}
\end{align*}
Thus the product becomes
\[
    \prod_{p \mid t} D_p(r) =
    \prod_{p \mid t} (p(1-p^2))^{-1} \prod_{\substack{p \mid t, \\ p \mid r}} (1-p^3)
    \prod_{\substack{p \mid t, \\ p \nmid r}} \left(1-
         \epsilon_p \left( 
         \frac{-ru_p}{p}
         \right)p^{3/2} 
         \right).
\]
Denoting by $\omega_3(x)$ the number of primes $p \mid x$ with $p \equiv 3 \bmod 4$, we have 
\begin{align*}
    \prod_{\substack{p \mid t, \\ p \nmid r}} \left(1-
         \epsilon_p \left( 
         \frac{-ru_p}{p}
         \right)p^{3/2} 
         \right) 
         &=
        \sum_{\substack{x \mid \tau,
          \\
           (x,r)=1
          }}
          \mu(x)
          \prod_{p \mid x}
          \left(
          \epsilon_p \left( 
         \frac{-ru_p}{p}
         \right)p^{3/2}  
         \right)
          \\
          &=
        \sum_{\substack{x \mid \tau,
          \\
           (x,r)=1
          }}
          \mu(x)
          x^{3/2}
          \left( \frac{-r}{x}
          \right)
           \prod_{p \mid x}
           \Bigg(
          \epsilon_p 
          \prod_{\substack{p' | kt, \\ p' \nmid x}}  
          \left(
          \frac{p'}{p}
          \right)
          \prod_{\substack{p' | x, \\p' \neq p}}  
          \left(
          \frac{p'}{p}
          \right)
         \Bigg)
          \\
         & =
          \sum_{\substack{x \mid \tau,
          \\
           (x,r)=1
          }}
          \mu(x) x^{3/2}
          \left( \frac{-r}{x}
          \right)
          \left( \frac{ \kappa \tau/x}{x}
          \right)
          (-1)^{\binom{\omega_3(x)}{2}}
          \prod_{p \mid x} \epsilon_p \\
          & =
          \sum_{\substack{x \mid \tau 
          }}
          \mu(x) x^{3/2}
          \left( \frac{-r \kappa \tau/x}{x}
          \right)\epsilon_x . 
\end{align*}
Here we used quadratic reciprocity,
\[
    \left(
    \frac{p}{p'}
    \right)
     \left(
    \frac{p'}{p}
    \right)
    = \begin{cases}
        -1 & \text{ if } p\equiv p' \equiv 3 \bmod 4, \\
        1 & \text{ otherwise},
    \end{cases}
\]
which implies that 
\[
    \prod_{p \mid x} \prod_{\substack{ p' \mid x \\ p' \neq p}} \left( \frac{p'}{p}
    \right) = (-1)^{\binom{\omega_3(x)}{2}},
\]
and that 
\[
(-1)^{\binom{\omega_3(x)}{2}} \prod_{p \mid x} \epsilon_p = (-1)^{\omega_3(x)(\omega_3(x)-1)/2}i^{\omega_3(x)} = i^{\omega_3(x)^2} = \epsilon_x. \qedhere
\]

\end{proof}

\subsection{Asymptotics for the contribution that is non-oscillatory after Poisson summation} \label{section:asymptotics_after_poisson}
 For any odd $k \ge 1$ that is not a perfect square and any positive integer $t$ with $ (2k,t)=1$, recall the definitions of  $k_1,k_2,\kappa_1,\kappa_2,\kappa$ and $\tau$  as well as $A_{k,t}, B_{k,t}, C_{k_2}$ and $D_{k,t}$
 from Lemma~\ref{lem:evaluationqsum}. 
 Let us choose an integer $r \in \Z$ of the form 
 \begin{equation}\label{eq:rforsecondorderterm}
         r = \sigma 2^j \ell^2 \kappa_2 \frac{\tau}{x_0}, 
 \end{equation}
 where $\sigma \in \{ \pm 1\}$, $j \ge 0$, $x_0$ divides $\tau$ and $\ell$ is an odd positive integer with $ (\ell, x_0 \kappa_1)=1$. These $r$ correspond to the non-oscillatory terms that contribute the secondary term for the asymptotic in Lemma~\ref{lem:mainlemmaboundedheight}.

    The goal of this section is to analyse the contribution from $A_{k,t}(r) \lambda_{k,t} (r)$, where $\lambda_{k,t}(r)$ is the product of $B_{k,t}(r), \, C_{k_2}(r)$ and the $x = x_0$ summand from $D_{k,t}(r)$, i.e. 
    \begin{equation}\label{eq:def_lambda_kt}
            \lambda_{k,t} (r):=
     \frac{ \varphi(\kappa_2)}{\kappa_2} \prod_{p \mid \tau} \frac{(1-p^3)}{p(1-p^2)} \prod_{\substack{p \mid x_0 }}\! \frac{(-1)}{1-p^3} \,
    \frac{ x_0^{3/2} }{\kappa_1^{1/2} }
    \epsilon_{x_0 \kappa_1}
           \left(\frac{- \sigma 2^j  }{\kappa_1 x_0}
    \right),
    \end{equation}
    with $r \in \mathcal{R}_{k,t}$, where we set 
    \[
    \mathcal{R}_{k,t} =
\left\{ r:~ r = \sigma 2^j \ell^2 \kappa_2 \frac{\tau}{x_0} \text{ is of the form \eqref{eq:rforsecondorderterm}} \right\}.
\]
We will establish the following asymptotic formula.

\begin{lemma} \label{lem:mainsecondordertermstatement}
   Let $Q_{k,t}$ be defined as in Lemma~\ref{lem:possionsummationapplication}.
    Then
    \begin{align*}
    X \sum_{\substack{k\le Y, \\ k \neq \square, \\ \text{$k$ odd}}}  \frac{d(k)}{k}   
    \sum_{\substack{t \le \sqrt{X}, \\  (t,2k)=1}}  \frac{2^{\omega(t)}}{t^3}  \sum_{r \in \mathcal{R}_{k,t}}  \widehat{w}(\tfrac{Xr}{Q_{k,t}t^2}) A_{k,t}(r) \lambda_{k,t}(r) = \frac{-48 \log(2)}{\pi \zeta(3)} \sqrt{X} + O\big( E(X,Y,Z)
      \big),
    \end{align*}
    where 
    \[
    E(X,Y,Z) \ll_\varepsilon (XYZ)^\varepsilon (X Y^{-1/2} + Z + (X/Z)^{1/2}).
    \]
\end{lemma}

For the proof of Lemma~\ref{lem:mainsecondordertermstatement}, we need some preparation. We will use Mellin inversion to evaluate sums over $\widehat{w}$. For $z \in \C$, we set
\begin{equation*}
     \mathcal{W}(z) = \int_0^1 w(u) u^{z-1} \mathrm{d}u.
\end{equation*}
Let $\sigma>0$ be a positive number. 
 Using integration by parts as in the proof of \eqref{eq:fouriertrafobound} yields, after a change of variables $u = e^{-t}$,
\begin{equation*}
        |\mathcal{W}(\sigma+i\tau)|
    \ll_\sigma
    \min\left\{1,\frac{1}{|\tau|},\frac{Z}{|\tau|^2}\right\}.
\end{equation*}
Consequently, for $0\le \alpha<1$, we have
\begin{equation}\label{eq:boundW}
    \int_{-\infty}^{\infty}
    (1+|\tau|)^\alpha |\mathcal{W}(\sigma+i\tau)|\,\mathrm{d}\tau
    \ll_{\sigma,\alpha}
    \begin{cases}
        \log Z & \text{ if } \alpha=0,\\
        Z^\alpha & \text{ if } 0<\alpha<1.
    \end{cases}
\end{equation}
Since $0\le w\le1$ and
$w=1$ on $[Z^{-1},1-Z^{-1}]$, we have
\begin{equation}\label{eq:w1/2}
    \mathcal{W}(1/2) = 2 + O(Z^{-1/2}).
\end{equation}

   Moreover, we record the following identity.
\begin{lemma}\label{lem:mellin}
    Let $c \in (1/2,1)$ be a real number, let $v > 0$ and let $m$ be an odd integer. Then we have 
    \begin{equation}\label{eq:what}
        \sum_{\substack{\ell \ge 1, \\
        (\ell,m)=1}} \widehat{w}(v \ell^2)
        = 
        \frac{1}{2 \pi i} 
        \int_{c-i\infty}^{c+i \infty} \frac{\Gamma(s) \zeta(2s)}{(2\pi)^s}
        \prod_{p \mid m} (1-p^{-2s})\,
        e^{-\frac{\pi i s}{2} } \mathcal{W}(1-s)v^{-s}
        \mathrm{d}s,
    \end{equation}
    and the integral is absolutely convergent. 
\end{lemma}

\begin{proof}
   Since all derivatives of  $\widehat{w}$ are rapidly decreasing by \eqref{eq:whatderivativebound},  Mellin inversion (see \cite[\S5.1]{multntmontvaugh}) is applicable and gives
    \[
    \widehat{w}(v\ell^2) = 
    \frac{1}{2 \pi i} \int_{c-i\infty}^{c+i \infty} \bigg(
    \int_0^\infty \widehat{w}(u) u^{s-1}\mathrm{d}u
    \bigg) (v\ell^2)^{-s} \mathrm{d}s.
    \]
   By the formula for the oscillatory Gamma integral (see, for example,
\cite[Appendix~C]{multntmontvaugh}), for $0<\Re(s)<1$ and $x>0$, we have
\[
\int_0^\infty e^{-2\pi i x u}u^{s-1}\,du
=
\Gamma(s)(2\pi x)^{-s}e^{-\pi i s/2}.
\]
Therefore 
\[
\int_0^\infty \widehat{w}(u) u^{s-1}\mathrm{d}u = \frac{\Gamma(s)}{(2\pi)^{s}} e^{-\pi i s/2} \mathcal{W}(1-s) 
\]

For $\Re(s)>1/2$, the sum 
\[
    \sum_{\substack{\ell \ge 1, \\ (\ell,m)=1}} \ell^{-2s} = \zeta(2s)  \prod_{p \mid m} (1-p^{-2s})
\]
converges absolutely and so substituting this identity gives the claimed formula \eqref{eq:what}.
By Stirling's formula (see, for example, \cite[Appendix~C]{multntmontvaugh}),
for $s=c+i\tau$ we have
\[
    \left|\Gamma(s)e^{-\pi i s/2}\right|
    \ll_c (1+|\tau|)^{c-1/2}.
\]
This, in combination with \eqref{eq:boundW} and the bound $|\zeta(2s)| \le \zeta(2c)<\infty$, shows the absolute convergence of the integral. 
\end{proof}

Lastly, let  $\chi$ denote any character modulo 8 and let $t$ denote any positive odd integer. Define the auxiliary function 
\begin{align}\label{eq:defofFchits}
F_{\chi,t,Y}(s) =
     \sum_{\substack{k\le Y, \\ k \neq \square, \\ (k,2t)=1}} & \frac{d(k)}{k} 
    \frac{\varphi(\kappa_2)}{\kappa_2}
\kappa_1^{s-\frac{1}{2}}\chi(\kappa_1) \prod_{p\mid\kappa_1}(1-p^{-2s}).
\end{align}
Then we have the following result.
\begin{lemma}\label{lem:nonprincipalbound}
    Fix $\delta \in (0,\frac{1}{4})$ and let $\frac{1}{2}+\delta \le \sigma \le 1-\delta$ and $\chi$ be a nonprincipal character modulo 8. Then 
    $
       F_{\chi,t,Y}(\sigma + i \tau) \ll_{\delta,\varepsilon} t^\varepsilon (1 + |\tau|)^{ \frac{1}{2}+\varepsilon}.
    $
\end{lemma}
\begin{proof}
We begin with simplifying the sum. The restriction $k \neq \square$ only gives a contribution of
$
\sum_{m\ge1}{d(m^2)}{m^{-2}}\ll1.
$
Since $\chi$ is a quadratic character,
we have $\chi(\kappa_1)=\chi(k)$.
The function
\[
a_s(k):=
k^{1/2-s}
\frac{\varphi(\kappa_2)}{\kappa_2}
\kappa_1^{s-1/2}
\prod_{p\mid\kappa_1}(1-p^{-2s})
\mathbf{1}_{(k,2t)=1}
\]
is multiplicative. Setting
$
A_s(m):=\sum_{d\mid m}\mu(d)a_s(m/d),
$
we have by
M\"obius inversion
\begin{align*}
F_{\chi,t,Y}(s)
=
\sum_{k\le Y}
\frac{d(k)\chi(k)}{k^{3/2-s}}
\sum_{m\mid k}A_s(m)+O(1) =
\sum_{m\le Y}
\frac{A_s(m)\chi(m)}{m^{3/2-s}}
\sum_{\ell\le Y/m}
\frac{d(m\ell)\chi(\ell)}
{\ell^{3/2-s}}
+O(1).
\end{align*}
We extend the definition of $\kappa_1, \kappa_2$ given in  Lemma~\ref{lem:evaluationqsum} to all odd integers $k$, using the same formulas.
Using the identity
\[
d(m\ell)=
\sum_{h\mid(m,\ell)}
\mu(h)d(m/h)d(\ell/h),
\]
and writing $\ell=hn$, we obtain
\begin{align*}
F_{\chi,t,Y}(s)=
\sum_{m\le Y}
\frac{A_s(m)\chi(m)}{m^{3/2-s}}
\sum_{h\mid m}
\frac{\mu(h)d(m/h)\chi(h)}
{h^{3/2-s}} 
\sum_{n\le Y/(mh)}
\frac{d(n)\chi(n)}{n^{3/2-s}}
+O(1).
\end{align*}

We now bound the innermost sum. By the
work of
Fujii--Gallagher--Montgomery~\cite[Theorem~1]{fujii_gallagher_montgomery}, we have
\[
\sup_{u\ge1}
\Bigg|
\sum_{n\le u}\chi(n)n^{i\tau}
\Bigg|
\ll (1+|\tau|)^{1/2}\log(2+|\tau|).
\]
Dirichlet's hyperbola method therefore gives
\[
\sum_{n\le u}d(n)\chi(n)n^{i\tau}
\ll
u^{1/2}(1+|\tau|)^{1/2}
\log^2(2+|\tau|).
\]
Since $\frac{3}{2}-\sigma\ge\frac{1}{2}+\delta$,
partial summation yields
\[
\sup_{M\ge1}
\left|
\sum_{n\le M}
\frac{d(n)\chi(n)}{n^{3/2-s}}
\right|
\ll_{\delta,\varepsilon}
(1+|\tau|)^{1/2+\varepsilon}.
\]
It follows that
\begin{align*}
F_{\chi,t,Y}(s)
\ll_{\delta,\varepsilon}\ &
(2+|\tau|)^{1/2+\varepsilon}
\sum_{m\le Y}
\frac{|A_s(m)|}{m^{3/2-\sigma}}
\sum_{h\mid m}
\frac{|\mu(h)|d(m/h)}
{h^{3/2-\sigma}}
+1.
\end{align*}

It remains to bound the sum over $m$.
Since $a_s$ is multiplicative, so is $A_s$.
For $p\nmid2t$, we have
\[
a_s(p^r)=
\begin{cases}
p^{(r-1)(1/2-s)}(1-p^{-2s})
&\text{if }r\text{ is odd},\\
p^{r(1/2-s)}(1-p^{-1})
&\text{if }r\text{ is even}.
\end{cases}
\]
Thus
\[
A_s(p)=a_s(p)-1=-p^{-2s},
\]
while for $r\ge2$ we have
\[
A_s(p^r)=a_s(p^r)-a_s(p^{r-1})
\ll1,
\]
uniformly for $\sigma\ge1/2+\delta$.
For $p\mid2t$, we instead have
$A_s(p)=-1$ and $A_s(p^r)=0$ for $r\ge2$.

Therefore, by multiplicativity,
\begin{align*}
&\sum_{m\ge1}
\frac{|A_s(m)|}{m^{3/2-\sigma}}
\sum_{h\mid m}
\frac{|\mu(h)|d(m/h)}
{h^{3/2-\sigma}}\\
&\quad\ll_\delta
\prod_{p\nmid2t}
\left(
1+O_\delta\left(
p^{-3/2-\sigma}+p^{-3+2\sigma}
\right)
\right)
\prod_{p\mid2t}
\left(1+p^{\sigma-3/2}\right)^2.
\end{align*}
Since $\sigma\le1-\delta$, the first
product is bounded uniformly, while
the second is
$\ll_{\delta,\varepsilon}t^\varepsilon$.
Hence
\[
\sum_{m\ge1}
\frac{|A_s(m)|}{m^{3/2-\sigma}}
\sum_{h\mid m}
\frac{|\mu(h)|d(m/h)}
{h^{3/2-\sigma}}
\ll_{\delta,\varepsilon}t^\varepsilon\qedhere
\]
\end{proof}

\begin{proof}[Proof of Lemma~\ref{lem:mainsecondordertermstatement}]
By the parametrisation
\[
    r=\sigma 2^j\ell^2\kappa_2\frac{\tau}{x_0},
\]
and since $Q_{k,t}= 16 \kappa_1 \kappa_2 \tau $, we have
\begin{align*}
\sum_{r\in\mathcal{R}_{k,t}} \!\!\!
\widehat{w} \left(\tfrac{Xr}{Q_{k,t}t^2} \right) \!
A_{k,t}(r) \lambda_{k,t}(r) 
 =
\sum_{x_0\mid\tau}
\sum_{\sigma\in\{\pm1\}}
\sum_{j\ge0}
A_{k,t}
\Big(\sigma2^j\kappa_2\frac{\tau}{x_0}\Big)
\lambda_{k,t}
\Big(\sigma2^j\kappa_2\frac{\tau}{x_0}
\Big) \!\!\!\!
\sum_{\substack{\ell\ge1,\\ (\ell,2x_0\kappa_1)=1}} \!\!\!\!\!
\widehat{w}\left(
\tfrac{X\sigma2^j\ell^2}{16\kappa_1t^2x_0}
\right).
\end{align*}
Here we used the fact that both $A_{k,t}(r)$ and $\lambda_{k,t}(r)$ are independent
of $\ell$.
Setting
\[
\Sigma_{k,t}(\sigma,j,x_0)
:=
\sum_{\substack{\ell\ge1, \\ (\ell,x_0\kappa_1)=1}}
\widehat{w}\left(
\tfrac{X\sigma2^j\ell^2}{16\kappa_1t^2x_0}
\right),
\]
we have
\[
\sum_{\substack{\ell\ge1,\\ (\ell,2x_0\kappa_1)=1}}
\widehat{w}\left(
\tfrac{X\sigma2^j\ell^2}{16\kappa_1t^2x_0}
\right)
=
\Sigma_{k,t}(\sigma,j,x_0)-\Sigma_{k,t}(\sigma,j+2,x_0).
\]
Since $x_0\kappa_1$ is odd, for $j\ge2$ we have
\[
\left(\frac{-\sigma2^j}{x_0\kappa_1}\right)
=
\left(\frac{-\sigma2^{j-2}}{x_0\kappa_1}\right),
\]
and therefore
\[
\lambda_{k,t} \Big(
\sigma2^j\kappa_2\frac{\tau}{x_0}
\Big)
=
\lambda_{k,t}\Big(
\sigma2^{j-2}\kappa_2\frac{\tau}{x_0}
\Big).
\]
Hence, after reindexing the second term, we get 
\begin{align*}
&\sum_{r\in\mathcal{R}_{k,t}}
\widehat{w} \left(\tfrac{Xr}{Q_{k,t}t^2} \right)
A_{k,t}(r) \lambda_{k,t}(r) \\ 
&\qquad \quad = \sum_{
\substack{
x_0\mid\tau, \\
j\ge0, \\
\sigma\in\{\pm1\}
}
}
\left(
A_{k,t}\Big(\sigma2^j\kappa_2\frac{\tau}{x_0}\Big)
-
\mathbf{1}_{j\ge2}
A_{k,t}\Big(\sigma2^{j-2}\kappa_2\frac{\tau}{x_0}\Big)
\right)
\lambda_{k,t}\Big(
\sigma2^j\kappa_2\frac{\tau}{x_0}
\Big)
\Sigma_{k,t}(\sigma,j,x_0).
\end{align*}

For $j\ge6$, the expression in parentheses vanishes. Indeed, in
this case both $2^j$ and $2^{j-2}$ are divisible by $16$ and hence
\[
A_{k,t}\Big(
\sigma2^j \kappa_2 \frac{\tau}{x_0}
\Big)
=
A_{k,t}\Big(
\sigma2^{j-2}\kappa_2\frac{\tau}{x_0}
\Big)
=
\frac{70}{\pi^2}.
\]
It follows that only the terms with $0\le j\le5$ remain.
For odd $m$ and $0\le j\le5$, we define
\begin{equation}\label{def:gamma_j}
\gamma_j(m)
:=
\frac{1}{\pi^2 8}
\sum_{a\bmod16} b_a
\left(
e_{16}(2^jam)
-
\mathbf{1}_{j\ge2}e_{16}(2^{j-2}am)
\right),
\end{equation}
so that
\[
A_{k,t}\Big(
2^j\kappa_2\frac{\tau}{x_0}
\Big)
-
\mathbf{1}_{j\ge2}
A_{k,t}\Big(
2^{j-2}\kappa_2\frac{\tau}{x_0}
\Big)
=
\frac{1}{2}\gamma_j(x_0\kappa_1).
\]
Finally, we have
\[
A_{k,t}(-r)=\overline{A_{k,t}(r)},
\qquad
\lambda_{k,t}(-r)=\overline{\lambda_{k,t}(r)},
\qquad
\Sigma_{k,t}(-1,j,x_0)
=
\overline{\Sigma_{k,t}(1,j,x_0)}.
\]
Thus the terms corresponding to $\sigma=\pm1$ combine to give
\begin{align*}
\sum_{r\in\mathcal{R}_{k,t}}
\widehat{w}\left(\tfrac{Xr}{Q_{k,t}t^2}\right)
A_{k,t}(r)\lambda_{k,t}(r)
=
\Re 
\sum_{x_0\mid\tau}
\sum_{j=0}^{5}
\gamma_j(x_0\kappa_1)
\lambda_{k,t}\Big(
2^j\kappa_2\frac{\tau}{x_0}
\Big)
\Sigma_{k,t}(j,x_0)
,
\end{align*}
where we abbreviate $\Sigma_{k,t}(j,x_0) = \Sigma_{k,t}(1,j,x_0) $.

We now evaluate $\Sigma_{k,t}(j,x_0)$ using Mellin inversion identities. Indeed, by Lemma~\ref{lem:mellin}, we have for any $c \in (1/2,1)$ that 
\begin{align*}
\Sigma_{k,t}(j,x_0)
&=
\frac{1}{2\pi i}
\int_{c-i\infty}^{c+i\infty}
\frac{\Gamma(s)\zeta(2s)}{(2\pi)^s}
\prod_{p\mid x_0\kappa_1}(1-p^{-2s})
e^{-\pi i s/2}\mathcal{W}(1-s) 
\frac{(\kappa_1t^2x_0)^s}{(X2^{j-4})^s}
\mathrm{d}s .
\end{align*}
We now observe that all dependence on $j$ in the sum over $\gamma_j(x_0\kappa_1)
\lambda_{k,t}\Big(
2^j\kappa_2\frac{\tau}{x_0}
\Big)
\Sigma_{k,t}(j,x_0)$  comes  from $\gamma_j(x_0\kappa_1)$ as well as contributions from the integral representation of $\Sigma_{k,t}(j,x_0)$ in the form of $2^{-(j-4)s}$ and from 
$\lambda_{k,t}\Big(
2^j\kappa_2\frac{\tau}{x_0}
\Big)$ in the form (see \eqref{eq:def_lambda_kt}) of
\[
  \epsilon_{x_0 \kappa_1}
           \left(\frac{- 2^j  }{\kappa_1 x_0}
    \right).
\]
Exchanging the order of integration and summation over $j$ thus yields the coefficient 
\[
    \sum_{j=0}^5
    \epsilon_{x_0 \kappa_1}
           \left(\frac{- 2^j  }{\kappa_1 x_0}
    \right) 
    \gamma_j(x_0 \kappa_1) 2^{-(j-4)s},
\]
which only depends on the value of the odd integer $x_0 \kappa_1$ modulo 8.

Denoting by $\mathcal{C} = \{\chi_0, \chi_{-4}, \chi_8, \chi_{-8} \}$ the four real Dirichlet characters modulo 8, with $\chi_0$ principal, we set 
\begin{equation}\label{eq:def_c_chi}
    c_{\chi}(s) = \frac{1}{4}\sum_{a \in (\Z / 8\Z)^{\times}} \chi(a)
     \sum_{j=0}^5
    \epsilon_{a}
           \left(\frac{- 2^j  }{a}
    \right) 
    \gamma_j(a) 2^{-(j-4)s}
\end{equation}
for $\chi \in \mathcal{C}$.
In particular, we have 
\[
    c_{\chi_0} (s) = \frac{1-i}{\pi^2} 
      \left(
        -3\cdot 2^{4s+1}
        -2^{3s+9/2}
        -5\cdot 2^{2s+2}
        + 2^{s+9/2}
        +96
    \right),
\]
which allows us to record 
\begin{align}\label{eq:catprincipal}
    c_{\chi_0} \left(\tfrac{1}{2}\right)=0, 
    \qquad
    c_{\chi_0}'\left(\tfrac{1}{2}\right) = \frac{(-1+i)336 \log 2}{\pi^2}.
\end{align}

Fourier inversion yields
\[
    \sum_{j=0}^5
    \epsilon_{x_0 \kappa_1}
           \left(\frac{- 2^j  }{\kappa_1 x_0}
    \right) 
    \gamma_j(x_0 \kappa_1) 2^{-(j-4)s} = \sum_{\chi \in \mathcal{C}} c_\chi(s) \chi(x_0 \kappa_1).
\]
Employing the definition of $\lambda_{k,t}$ in \eqref{eq:def_lambda_kt}, we thus arrive at 
\begin{align*}
\sum_{r\in\mathcal{R}_{k,t}}
\widehat{w} 
\left(\tfrac{Xr}{Q_{k,t}t^2}\right) &
A_{k,t}(r)\lambda_{k,t}(r) 
\\
=
\Re 
\sum_{\chi \in \mathcal{C}} &
\sum_{x_0\mid\tau}\chi(x_0)\chi(\kappa_1)
  \frac{ \varphi(\kappa_2)}{\kappa_2} \prod_{p \mid \tau} \frac{(1-p^3)}{p(1-p^2)} \prod_{\substack{p \mid x_0 }}\! \frac{(-1)}{1-p^3} \,
    \frac{ x_0^{3/2} }{\kappa_1^{1/2} } \\
& \qquad \cdot \frac{1}{2\pi i}
\int_{c-i\infty}^{c+i\infty}
\frac{\Gamma(s)\zeta(2s)}{(2\pi)^s}
\prod_{p\mid x_0\kappa_1}(1-p^{-2s})
e^{-\pi i s/2} \mathcal{W}(1-s) 
\frac{(\kappa_1t^2x_0)^s}{X^s} c_\chi(s)
\mathrm{d}s 
.
\end{align*}

% Recalling that the ultimate goal of this lemma is to give asymptotics for 
% \[
%   X \sum_{\substack{k\le Y, \\ k \neq \square, \\ \text{$k$ odd}}}  \frac{d(k)}{k}   
%     \sum_{\substack{t \le \sqrt{X}, \\  (t,2k)=1}}  \frac{2^{\omega(t)}}{t^3}  \sum_{r \in \mathcal{R}_{k,t}}  \widehat{w}(\tfrac{Xr}{Q_{k,t}t^2}) A_{k,t}(r) \lambda_{k,t}(r),
% \]
% 
We now sum over $t$. Since the involved sums are finite, we may exchange the order of summation and integration. For fixed $k$ and $\chi$, we simplify 
\begin{align*} 
    \sum_{\substack{t \le \sqrt{X}, \\  (t,2k)=1}}  \frac{2^{\omega(t)}}{t^{3-2s}} &
     \prod_{p \mid \tau} \frac{(1-p^3)}{p(1-p^2)} 
    \sum_{x_0\mid\tau} \chi(x_0)  x_0^{s+3/2} 
 \prod_{\substack{p \mid x_0 }}\! \frac{(-1)}{1-p^3} 
\prod_{p\mid x_0\kappa_1}(1-p^{-2s}) \\
&=\prod_{p\mid\kappa_1}(1-p^{-2s})
\sum_{\substack{t\le\sqrt X, \\ (t,2k)=1}}
\frac{2^{\omega(t)}}{t^{3-2s}}
\prod_{p\mid\tau}\frac{p^3-1}{p(p^2-1)}
\sum_{x_0\mid\tau}
\chi(x_0)x_0^{s+3/2}
\prod_{p\mid x_0}\frac{1-p^{-2s}}{p^3-1}.
\end{align*}
By multiplicativity, we get 
\[
    \sum_{x_0\mid\tau}
\chi(x_0)x_0^{s+3/2}
\prod_{p\mid x_0}\frac{1-p^{-2s}}{p^3-1} = \prod_{p\mid t}
\left(
1+\frac{\chi(p)p^{s+3/2}(1-p^{-2s})}{p^3-1}
\right).
\]
With $F_{\chi,t,Y}(s)$ given in \eqref{eq:defofFchits}, we define
\begin{align*}
D_{\chi,X,Y}(s) = 
\sum_{\substack{t \le \sqrt{X}, \\  t \text{ odd}}}  \frac{2^{\omega(t)}}{t^{3-2s}}
\prod_{p\mid t}
\frac{p^3-1}{p(p^2-1)}
\left(
1+
\frac{\chi(p)p^{s+3/2}(1-p^{-2s})}{p^3-1}
\right)
F_{\chi,t,Y}(s)
\end{align*}
which gives  
\begin{align*}
     X \sum_{\substack{k\le Y, \\ k \neq \square, \\ \text{$k$ odd}}}  & \frac{d(k)}{k}   
    \sum_{\substack{t \le \sqrt{X}, \\  (t,2k)=1}}  \frac{2^{\omega(t)}}{t^3}  \sum_{r \in \mathcal{R}_{k,t}}  \widehat{w}(\tfrac{Xr}{Q_{k,t}t^2}) A_{k,t}(r) \lambda_{k,t}(r) \\
    & \qquad \qquad = \Re
\sum_{\chi \in \mathcal{C}}
\frac{1}{2\pi i}
\int_{c-i\infty}^{c+i\infty}
\frac{\Gamma(s)\zeta(2s)}{(2\pi)^s}
e^{-\pi i s/2}
\mathcal{W}(1-s)
c_\chi(s)
 D_{\chi,X, Y}(s)
X^{1-s} \mathrm{d} s
\end{align*}
with $c \in \left( \frac{1}{2},1\right)$. 
For $\chi \neq \chi_0$ and $s = c + i \tau$, with $c=1-\varepsilon$, Lemma~\ref{lem:nonprincipalbound} and the divisor bound yield
\[
|D_{\chi,X,Y}(s)| \ll_\varepsilon \sum_{t \le \sqrt{X}} t^{-3 + 2 c + \varepsilon/2} (1 + |\tau|)^{1/2+\varepsilon/2} \ll_\varepsilon (1 + |\tau|)^{1/2+\varepsilon/2}.
\]
By Stirling's formula (see \cite[Appendix C]{multntmontvaugh}), we have on $\Re(s)=c$
\[
 \left|\Gamma(s)e^{-\pi is/2}\zeta(2s)\right| \ll_c (1+|\tau|)^{c-1/2}.
\]
Combining these bounds with \eqref{eq:boundW}, applied with $\alpha=c+\varepsilon/2<1$, the full contribution from $\chi \neq \chi_0$ can be bounded by
\begin{equation}\label{eq:contribution_nontrivial_character}
        \ll_\varepsilon X^{1-c} Z^{c+\varepsilon/2}
        \ll_\varepsilon X^{\varepsilon} Z^{1+\varepsilon}.
\end{equation}

It remains to consider $\chi = \chi_0$. 
Here, we first remove the condition $k \neq \square$, i.e. we replace $D_{\chi_0,X,Y}$ with 
\begin{align*}
\widetilde{D}_{\chi_0,X,Y}(s) = 
\sum_{\substack{t \le \sqrt{X}, \\  t \text{ odd}}}  \frac{2^{\omega(t)}}{t^{3-2s}}
\prod_{p\mid t}
\frac{p^3-1}{p(p^2-1)}
\left(
1+
\frac{\chi(p)p^{s+3/2}(1-p^{-2s})}{p^3-1}
\right)
\widetilde{F}_{\chi,t,Y}(s),
\end{align*}
where 
\[
\widetilde{F}_{\chi,t,Y}(s)
=
 \sum_{\substack{k\le Y,  \\ (k,2t)=1}}  \frac{d(k)}{k} 
    \frac{\varphi(\kappa_2)}{\kappa_2}
\kappa_1^{s-\frac{1}{2}}\chi(\kappa_1) \prod_{p\mid\kappa_1}(1-p^{-2s})
\]
As in the proof of Lemma~\ref{lem:nonprincipalbound}, we have $F_{\chi,t,Y}(s) - \widetilde{F}_{\chi,t,Y}(s) = O(1) $. For $\Re(s) = c = 1-\varepsilon$, this yields
\[
    \widetilde{D}_{\chi_0,X,Y}(s)
    -D_{\chi_0,X,Y}(s)
    \ll_\varepsilon
    \sum_{t\le\sqrt{X}}t^{-3+2c+\varepsilon}
    =
    \sum_{t\le\sqrt{X}}t^{-1-\varepsilon}
    \ll_\varepsilon 1.
\]
Using Stirling's formula and \eqref{eq:boundW} as before, this contributes an error term 
\begin{equation}\label{eq:contribtution_trivial_char_sq}
    \ll_\varepsilon X^{2\varepsilon} Z^{1+\varepsilon}.
\end{equation}

We shift the line of integration $\Re s = c$ to the left of the pole of $\zeta(2s)$ at $s=1/2$, to some fixed line $\Re s =a$, where $a \in (0,1/2)$. Since $\Gamma (s)$ has no poles for $\Re s>0$ and $c_\chi(s)$, $\mathcal{W}(1-s)$ and the finite sums defining $D_{\chi_0,X,Y}(s)$ are also holomorphic in this region, the simple pole of $\zeta(2s)$ is indeed the only singularity crossed. This pole is cancelled by the zero of $c_{\chi_0}(s)$, as determined in \eqref{eq:catprincipal}.
Therefore 
\begin{align*}
    \frac{1}{2\pi i}
\int_{c-i\infty}^{c+i\infty} &
\frac{\Gamma(s)\zeta(2s)}{(2\pi)^s}
e^{-\pi i s/2}
\mathcal{W}(1-s)
c_{\chi_0}(s)
 \widetilde{D}_{\chi_0,X, Y}(s)
X^{1-s} \mathrm{d} s  \\
&\qquad \qquad = \frac{1}{2\pi i}
\int_{a-i\infty}^{a+i\infty}
\frac{\Gamma(s)\zeta(2s)}{(2\pi)^s}
e^{-\pi i s/2}
\mathcal{W}(1-s)
c_{\chi_0}(s)
  \widetilde{D}_{\chi_0 ,X, Y}(s)
X^{1-s} \mathrm{d} s.
\end{align*}
 
For $\Re(s) = a < 1/2$, we may define the untruncated Dirichlet series $\widetilde{D} (s)$ as
\begin{align*}
     \widetilde{D} (s)  =  \sum_{\substack{k\ge 1,  \\ \text{$k$ odd}}}  \frac{d(k)}{k} 
    \frac{\varphi(\kappa_2)}{\kappa_2}
\kappa_1^{s-\frac{1}{2}} \prod_{p\mid\kappa_1}(1-p^{-2s}) \!\!
 \sum_{\substack{t \ge 1, \\  (t,2k)=1}}  \frac{2^{\omega(t)}}{t^{3-2s}}
\prod_{p\mid t}
\frac{p^3-1}{p(p^2-1)}
\left(
1+
\frac{p^{s+3/2}(1-p^{-2s})}{p^3-1}
\right),
\end{align*}
which converges absolutely. By the divisor bound $d(n)\ll_\varepsilon n^\varepsilon$, estimating the tails $k>Y$ and $t>\sqrt X$ gives
\[
     \widetilde{D}(s) =  \widetilde{D}_{\chi_0,X,Y}(s) + O_{a,\varepsilon}(Y^{-1/2+a+\varepsilon} + X^{-1 + a +\varepsilon}).
\]
By the convexity bound for the Riemann zeta function and Stirling's formula (see \cite[Appendix C]{multntmontvaugh}), we have on $\Re(s)=a$
\[
 \left|\Gamma(s)e^{-\pi is/2}\zeta(2s)\right| \ll_{a,\varepsilon} (1+|\Im s|)^\varepsilon. 
\]
Together with \eqref{eq:boundW}, this implies
\[
\int_{-\infty}^{\infty}  \left|\Gamma(a+iv)e^{-\pi i(a+iv)/2}\zeta(2(a+iv))  \mathcal{W}(1-a-iv)\right| \mathrm{d}v \ll_{a,\varepsilon} Z^\varepsilon. 
\]
Thus replacing $ \widetilde{D}_{\chi_0, X,Y}$ by $ \widetilde{D}$ in the integral over $\Re s = \varepsilon$ introduces an error term of order 
\begin{equation}\label{eq:contribution_replacecutoofbyfull}
    \ll_\varepsilon X^{1-\varepsilon}  Z^\varepsilon (Y^{-1/2+2\varepsilon} + X^{-1   +2\varepsilon}). 
\end{equation}

The series $\widetilde D(s)$ can be factored  into
\begin{equation}\label{eq:dtildeasproductwithlfunction}
\widetilde D(s)
 =  \zeta\left(\tfrac{3}{2}-s\right)^2 H(s),
\end{equation}
where
\[
H(s)= (1-2^{-3/2+s})^2 \prod_{p>2}
\frac{1-p^{2s-3}}{(1-p^{-2})^2}
\left(
(1+p^{-2})(1+p^{-3})
-2p^{s-7/2}-2p^{-s-3/2}
\right).
\]
The product converges locally uniformly
for $0<\Re s<1$, and hence $H$ is
holomorphic in this strip.
In particular, at $s=\frac{1}{2}$ the odd-prime factors
simplify to $1-p^{-3}$, giving
\begin{equation}\label{eq:Hat1/2}
    H\!\left(\tfrac{1}{2}\right)
    =
    \frac{1}{4}\prod_{p>2}(1-p^{-3})
    =
    \frac{2}{7\zeta(3)}.
\end{equation}

We now move the contour in the integral involving $\widetilde{D}(s)$ back to the line $\Re(s) = c = 1-\varepsilon$. The only pole crossed is the double pole at $s = \frac{1}{2}$, and thus  
    \begin{align*}
    \Re \frac{1}{2\pi i}
\int_{a-i\infty}^{a+i\infty} &
\frac{\Gamma(s)\zeta(2s)}{(2\pi)^s}
e^{-\pi i s/2}
\mathcal{W}(1-s)
c_{\chi_0}(s)
 \widetilde{D} (s)
X^{1-s} \mathrm{d} s  \\
&\qquad \qquad = \Re \frac{1}{2\pi i}
\int_{c-i\infty}^{c+i\infty}
\frac{\Gamma(s)\zeta(2s)}{(2\pi)^s}
e^{-\pi i s/2}
\mathcal{W}(1-s)
c_{\chi_0}(s)
  \widetilde{D} (s)
X^{1-s} \mathrm{d}s - \Re R,
\end{align*}
with $R = \operatorname{Res}_{s = 1/2} \left( 
\Gamma(s)\zeta(2s)(2\pi)^{-s}
e^{-\pi i s/2}
\mathcal{W}(1-s)
c_{\chi_0}(s)
  \widetilde{D} (s)
X^{1-s} 
\right).
$
By \eqref{eq:catprincipal},
\eqref{eq:dtildeasproductwithlfunction}
and \eqref{eq:Hat1/2}, we have
\[
\zeta(2s)c_{\chi_0}(s)\widetilde D(s)
=
-(1-i)\left(
\frac{48\log 2}{\pi^2\zeta(3)}
\left(s-\tfrac12\right)^{-2}
+
A\left(s-\tfrac12\right)^{-1}
+O(1)
\right)
\]
for some $A\in\mathbb R$. Indeed, $A$ is a real number since the functions
$c_{\chi_0}(s)/(1-i)$ and $H(s)$
only take real values for real $s$.
Since $$(1-i)e^{-\pi is/2}
=
-i\sqrt{2}
-\frac{\pi}{\sqrt{2}}
\left(s-\tfrac12\right)
+O\left((s-\tfrac12)^2\right)$$
and since $
\Gamma(s)(2\pi)^{-s}
\mathcal W(1-s)X^{1-s}
$
is real for real $s$, we obtain by  \eqref{eq:w1/2}
\begin{equation}\label{eq:valueofresidue}
\Re R
=
-\frac{48\log 2}{\pi^2\zeta(3)} \left(-\frac{\pi}{\sqrt{2}}\right) \frac{\Gamma(\tfrac12)}{(2\pi)^{1/2}} \mathcal W\left(\tfrac12\right)X^{1/2} =
\frac{48\log 2}{\pi\zeta(3)}X^{1/2}
+
O\left((X/Z)^{1/2}\right).
\end{equation}

It remains to bound the integral on $\Re s=c$.
By \eqref{eq:dtildeasproductwithlfunction},
the boundedness of $H(s)$ and $c_{\chi_0}(s)$
on this line, and the convexity bound for
the Riemann zeta function, we have
\[
\widetilde D(c+iv)
\ll_{c,\varepsilon}
(1+|v|)^{c-1/2+\varepsilon}.
\]
Moreover, Stirling's formula gives
\[
\Gamma(c+iv)e^{-\pi i(c+iv)/2}
\zeta(2c+2iv)
\ll_c (1+|v|)^{c-1/2}.
\]
Consequently, by \eqref{eq:boundW} and with $c = 1-\varepsilon$,
\begin{align*}
\int_{c-i\infty}^{c+i\infty}
\frac{\Gamma(s)\zeta(2s)}{(2\pi)^s}
e^{-\pi is/2}\mathcal W(1-s)
c_{\chi_0}(s)\widetilde D(s)X^{1-s}\,\mathrm{d}s
\ll_{c,\varepsilon}
X^{1-c}Z^{2c-1+\varepsilon}
\ll_\varepsilon (XYZ)^\varepsilon Z.
\end{align*}
Combining this with the error terms arising from \eqref{eq:contribution_nontrivial_character}, \eqref{eq:contribtution_trivial_char_sq}, \eqref{eq:contribution_replacecutoofbyfull} and \eqref{eq:valueofresidue}
completes the proof.
\end{proof}

\subsection{Bounds for the remaining contribution from the non-squares} \label{section:remaining_error_term_intheend}

It remains to bound 
\begin{equation}\label{eq:defE1}
E_1 :=
X \sum_{
     \substack{k \le Y, \\ k \neq \square, \\ k \text{ odd}
     }
     }
     \frac{d(k)}{k} \!\!\!
     \sum_{\substack{t \le \sqrt{X}, \\  (t,2k)=1}} \! \! \! \frac{2^{\omega(t)}}{t^3} 
         \sum_{\substack{r \in \Z , \\ 
         r \not \in \mathcal{R}_{k,t}}}
         \widehat{w}\big(\tfrac{X r}{Q_{k,t}t^2} \big)
         H_{k,t}(r),
\end{equation}
and 
\begin{equation}\label{eq:defE2}
E_2 :=
    X \sum_{
     \substack{k \le Y, \\ k \neq \square, \\ k \text{ odd}
     }
     }
     \frac{d(k)}{k} \!\!\!
     \sum_{\substack{t \le \sqrt{X}, \\  (t,2k)=1}} \! \! \! \frac{2^{\omega(t)}}{t^3} 
         \sum_{\substack{ 
         r  \in \mathcal{R}_{k,t}}}
         \widehat{w}\big(\tfrac{X r}{Q_{k,t}t^2} \big)
         (H_{k,t}(r) - A_{k,t}(r)\lambda_{k,t}(r))
\end{equation}
where $\mathcal{R}_{k,t}$ was defined in Lemma~\ref{lem:mainsecondordertermstatement}, $H_{k,t}(r)$ was defined in Lemma~\ref{lem:possionsummationapplication}, $A_{k,t}(r)$ was defined in Lemma~\ref{lem:evaluationqsum} and $\lambda_{k,t}(r)$ was defined in \eqref{eq:def_lambda_kt}. 

We use a consequence of Jutila's mean-square estimate
\cite[Theorem~1]{Jutila1975Mean}, which states that
    \begin{align} \label{eq:jutila2}
    \sum_{\substack{1\le |D| \le R, \\ D \text{ fundamental}, \\ D \neq 1}}
        \int_{-T}^T 
        \left|L 
        \left(
        \tfrac{1}{2} + iy, \left( 
        \tfrac{D}{\cdot}
        \right)
        \right)\right|^2 \mathrm{d}y \ll_\varepsilon 
        (RT)^{1+\varepsilon}
    \end{align}
    for any $R,T \ge 2$.
    Introducing some mild Gaussian smoothing, we have for any fundamental $D \neq 1$ and any small $\delta > 0$ that 
    \[
        \int_{-T}^T 
        \left|L 
        \left(
        \tfrac{1}{2} +\delta + iy, \left( 
        \tfrac{D}{\cdot}
        \right)
        \right)\right|^2  \mathrm{d}y  \ll 
        \int_{-\infty}^\infty 
        \left|L 
        \left(
        \tfrac{1}{2}+\delta + iy, \left( 
        \tfrac{D}{\cdot}
        \right)
        \right) e^{s^2/T^2} \right|^2   \mathrm{d}y,
    \]
    where we write $s = \tfrac{1}{2}+\delta + iy$.
     By the convexity principle for integral means on a vertical strip we get
    \begin{align*} 
        & \int_{-\infty}^\infty  
        \left|L 
        \left(
        \tfrac{1}{2} + \delta + iy, \left( 
        \tfrac{D}{\cdot}
        \right)
        \right)
         e^{s^2/T^2} \right|^2   \mathrm{d}y 
        \\ 
        &  \ \ 
        \ll 
        \left(\int_{-\infty}^\infty 
        \left|L 
        \left(
        \tfrac{1}{2} + iy, \left( 
        \tfrac{D}{\cdot}
        \right)
        \right)  \right|^2 e^{-2y^2/T^2}  \mathrm{d}y \right)^{1-\delta}
        \left(
        \int_{-\infty}^\infty 
        \left|L 
        \left(
        \tfrac{3}{2} + iy, \left( 
        \tfrac{D}{\cdot}
        \right)
        \right)  \right|^2 e^{-2y^2/T^2} \mathrm{d}y
        \right)^{\delta}.
    \end{align*}
Summing over $D$, applying H\"older's inequality with exponents $(1/(1-\delta), 1/\delta)$, splitting the range of integration into $[-T,T]$ and dyadic intervals $\pm [ 2^jT, 2^{j+1} T]$ and applying \eqref{eq:jutila2} yields  \begin{align*}
       \sum_{\substack{1\le |D| \le R, \\ D \text{ fundamental}, \\
       D \neq 1}} &
        \int_{-T}^T 
        \left|L 
        \left(
        \tfrac{1}{2} +\delta + iy, \left( 
        \tfrac{D}{\cdot}
        \right)
        \right)\right|^2  \mathrm{d}y  
        \\&\ll_{\varepsilon,\delta}
(RT)^\delta
\bigg(
(RT)^{1+\varepsilon}
+\sum_{j\ge0}e^{-4^j}
(2^jRT)^{1+\varepsilon}
\bigg)^{1-\delta}
\ll_{\varepsilon,\delta}(RT)^{1+\varepsilon}.
\end{align*}
    Setting $\delta = \varepsilon$ yields the following.
\begin{lemma}
\label{lem:lfunctionmeansquare}
    Let $T, R \ge 2$. Then 
     \[
        \sum_{\substack{1\le |D| \le R, \\ D \text{ fundamental}, \\ D \neq 1}}
        \int_{-T}^T 
        \left|L 
        \left(
        \tfrac{1}{2} +\varepsilon + iy, \left( 
        \tfrac{D}{\cdot}
        \right)
        \right)\right|^2  \mathrm{d}y  \ll_\varepsilon (RT)^{1+\varepsilon}.
    \]
\end{lemma}

We are now ready to bound the final error terms.

\begin{lemma} 
% \label{lem:remainingcancellation}
    With $E_1$ and $E_2$ defined as in \eqref{eq:defE1} and \eqref{eq:defE2}, we have 
    \[
|E_1| + |E_2| 
\ll_\varepsilon
         (XYZ)^{\varepsilon}Z.
\]
\end{lemma}

\begin{proof}
   For any $k,t$ in the summation range, let us recall the definitions of $k_1$, $k_2$, $\kappa$, $\kappa_1$, $\kappa_2$ and $\tau$ from Lemma~\ref{lem:evaluationqsum}, i.e. that we write $k=k_1 k_2^2$, 
    where $\tau = \rad(t)$, 
    \[
    k_1 = \prod_{2 \nmid \nu_p(k)} p^{\nu_p(k)}, \qquad k_2 = \prod_{2 \mid \nu_p(k)} p^{\nu_p(k)/2},
    \]
    and $\kappa_1 = \rad(k_1)$, $\kappa_2 = \rad(k_2)$.  Moreover, choosing a divisor $x_0 \mid \tau$, we consider all $r \in \Z$ that are \emph{not} of the form 
 \begin{equation*}
          \sigma 2^j \ell^2 \kappa_2 \frac{\tau}{x_0}, 
 \end{equation*}
 where $\sigma \in \{ \pm 1\}$, $j \ge 0$, $x_0$ divides $\tau$ and $\ell$ is an odd positive integer with $ (\ell, x_0 \kappa_1)=1$.
 Then, by decomposition $H_{k,t}(r) = A_{k,t}(r) B_{k,t}(r) C_{k_2}(r) D_{k,t}(r)$ provided by Lemma~\ref{lem:evaluationqsum}, we observe that to prove the lemma statement, it suffices to bound the contribution from the product of $A_{k,t}(r)$, $B_{k,t}$, $C_{k_2}(r)$ and the $x = x_0$ summand in \[
    D_{k,t}(r) = 
          \prod_{p \mid t} (p(1-p^2))^{-1} \prod_{\substack{p \mid t, \\ p \mid r}} (1-p^3)
          \sum_{\substack{x \mid \tau
          }}
          \mu(x) x^{3/2}
          \left( \frac{-r \kappa \tau/x}{x}
          \right)\epsilon_x,
\]
i.e. the contribution from 
\begin{align*}
\widetilde{\lambda}_{x_0}(r)=
    A_{k,t}(r)B_{k,t}(r)C_{k_2}(r)\prod_{p \mid t} (p(1-p^2))^{-1} \prod_{\substack{p \mid t, \\ p \mid r}} (1-p^3)
          \mu(x_0) x_0^{3/2}
          \left( \frac{-r \kappa \tau/x_0}{x_0}
          \right)\epsilon_{x_0}.
\end{align*}
 Since any non-zero integer $r$ can uniquely be represented as 
 \[
    r = \sigma 2^j q m^2,
 \]
 where $\sigma \in \{ \pm 1\}$, $j \ge 0$, $\mu(2q)^2 =1$ and $m$ is an odd positive integer, we have $(r,x_0\kappa_1)>1$, $r=0$ or $q \neq \kappa_2 \tau / x_0$. Using the rapid decay of $\widehat w$ from \eqref{eq:fouriertrafobound}, we may restrict the $r$-sum to $|r|\le (XYZ)^2$; the omitted tail is $\ll_A (XYZ)^{-A}$ for every $A>0$.

 If $ (r,x_0\kappa_1)>1$ or $r=0$, then $\widetilde{\lambda}_{x_0}(r)$ vanishes. Indeed, we have
 \[
    B_{k,t}(r) = \left(\frac{-r \kappa_2 \tau}{\kappa_1 }
    \right)
     \frac{\epsilon_{\kappa_1}}{\sqrt{\kappa_1}},
 \]
and the Jacobi symbol clearly vanishes for $ (r,\kappa_1)>1$. Similarly, if $ (r, x_0) > 1$, then the summand corresponding to $x = x_0$ in $D_{k,t}(r)$ is a multiple of 
$ \left( \frac{-r \kappa \tau/x}{x}
          \right)$ and thus
vanishes. 

Let us from now on assume that $(r, x_0 \kappa_1) = 1$ and consequently $q \neq \kappa_2 \tau / x_0$. We write $k_1 = \kappa_1 s_1^2$, so that $k = \kappa_1 s_1^2 k_2^2$ with $\rad(s_1) \mid \kappa_1$ and $(k_1,k_2)=(\kappa_1 s_1^2, k_2)=1$. Then, since $d$ is a multiplicative function, we have 
$d(k) = d(\kappa_1 s_1^2) d(k_2^2)$ and 
\[
    d(\kappa_1 s_1^2) = \prod_{p \mid \kappa_1} \left(
    2 + 2\nu_p(s_1)
    \right) 
    = 2^{\omega(\kappa_1)} \prod_{p \mid \kappa_1} \left(
    1 + \nu_p(s_1)
    \right) 
    = 
    2^{\omega(\kappa_1)} \prod_{p \mid s_1} \left(
    1 + \nu_p(s_1)
    \right).
\]

We observe that $\widetilde{\lambda}_{x_0}(r)$ is the product of 
\[
    \prod_{p \mid t} (p(1-p^2))^{-1} 
          \mu(x_0) x_0^{3/2}
          \epsilon_{x_0}\prod_{\substack{p \mid t, \\ p \mid r}} (1-p^3),
\]
which depends only on $r$, $\tau$ and $x_0$, with
\[
A_{k,t}(r)
    B_{k,t}(r)C_{k_2}(r) \left( \frac{-r \kappa_1 \kappa_2 \tau/x_0}{x_0}
          \right) .
\]
 By quadratic reciprocity, one can show  
\begin{align*}
     B_{k,t}(r) \left( \frac{-r \kappa_1 \kappa_2 \tau/x_0}{x_0}
          \right) &= \left(\frac{-r \kappa_2 \tau}{\kappa_1 }
    \right)
     \frac{\epsilon_{\kappa_1}}{\sqrt{\kappa_1}} \left( \frac{-r \kappa_1 \kappa_2 \tau/x_0 }{x_0}
          \right)  \\
          &= \frac{\epsilon_{\kappa_1x_0}}
{\epsilon_{x_0}\sqrt{\kappa_1}}
\left(
\frac{-r\kappa_2\tau/x_0}{\kappa_1}
\right)
\left(
\frac{-r\kappa_2\tau/x_0}{x_0}
\right).
\end{align*}

Fixing all other variables and considering only the sum over $\kappa_1$, we thus consider 
\begin{align}\label{eq:kappa1sumtobound}
    \sum_{\substack{1<\kappa_1 \le Y/(k_2s_1)^2, \\
    (\kappa_1, 2k_2t)=1, \\
    \rad(s_1)\mid \kappa_1}}
    \! \frac{2^{\omega(\kappa_1)}\mu(\kappa_1)^2}{\kappa_1^{3/2}} 
\left(
\frac{-r\kappa_2\tau/x_0}{\kappa_1}
\right)
\widehat{w}\big(\tfrac{X r}{16\tau \kappa_1\kappa_2 t^2} \big) \chi (\kappa_1),
\end{align}
where $\chi$ is a Dirichlet character modulo 8.
Here we used the fact that $A_{k,t}(r)$ and $\epsilon_{\kappa_1 x_0}$ only depend on $\kappa_1 \bmod 8$ and can therefore be expressed as a finite linear combination of characters modulo $8$. Indeed, for $\epsilon_{\kappa_1 x_0}$ this is immediate, while for $A_{k,t}(r)$ this follows from the fact that $b_{9N} = b_N$, where $b_N$ is defined in Theorem~\ref{thm:2-4}.

We use Perron's formula to estimate, uniformly for $1\le y\le Y/(k_2s_1)^2$, the sum 
\[
    \sum_{\substack{1<\kappa_1 \le y, \\
    (\kappa_1, 2k_2t)=1, \\
    \rad(s_1)\mid \kappa_1}}
    \! 2^{\omega(\kappa_1)}\mu(\kappa_1)^2
\left(
\frac{-r\kappa_2\tau/x_0}{\kappa_1}
\right) \chi(\kappa_1).
\] 
This sum vanishes unless $y \ge \max \{ 3, \rad (s_1)\}$ and $(\rad(s_1),2k_2t)=1$, so we may assume that these conditions hold from now on without loss of generality.
via partial summation, this will allow us to bound \eqref{eq:kappa1sumtobound}. Let us set 
\[
\eta(n) = \left(
\frac{-r\kappa_2\tau/x_0}{n}
\right) \chi(n),
\]
which, as a product of real quadratic characters, is a (possibly imprimitive) real quadratic character; in the range under consideration, the odd squarefree part $q$ of $r$ satisfies  
$
q\neq\kappa_2\frac{\tau}{x_0},$
so the odd squarefree part of
$-r\kappa_2\frac{\tau}{x_0}$
is nontrivial and consequently $\eta$ is nonprincipal. Moreover, we let 
\begin{align*}
    F(s) = \sum_{\substack{\kappa_1 \ge 1, \\
    (\kappa_1, 2k_2t)=1, \\
    \rad(s_1)\mid \kappa_1}}
    \! \frac{2^{\omega(\kappa_1)}\mu(\kappa_1)^2}{\kappa_1^s}
\eta(\kappa_1),
\end{align*}
which admits the factorisation 
\begin{align} \label{eq:facorisationF}
    F(s) =  2^{\omega(s_1)} \frac{\eta(\rad(s_1))}{\rad(s_1)^s} \prod_{p \nmid 2 s_1 k_2 t} \left(
1+ 2\eta(p)p^{-s} 
\right) =  2^{\omega(s_1)} \frac{\eta(\rad(s_1))}{\rad(s_1)^s} \frac{L(s, \eta)^2}{L(2s, \eta^2)} G(s), 
\end{align}
where 
\[
    G(s) = 
\prod_{p\nmid 2s_1k_2t} \left(
\left(1+2\eta(p)p^{-s}\right)
\frac{1-\eta(p)p^{-s}}
     {1+\eta(p)p^{-s}}
     \right)
\prod_{p\mid 2s_1k_2t}
\left(
\frac{1-\eta(p)p^{-s}}
     {1+\eta(p)p^{-s}}
     \right).
\]
Since 
\[
    \left(1+2\eta(p)p^{-s}\right)
\frac{1-\eta(p)p^{-s}}
     {1+\eta(p)p^{-s}} = 1+ O(p^{-2\Re(s)}),
\]
the infinite product defining $G(s)$ converges absolutely for $\Re(s)>1/2$. 

By Perron's formula \cite[Theorem~5.2 and Corollary~5.3]{multntmontvaugh} we have 
\begin{align*}
    \sum_{\substack{1<\kappa_1 \le y, \\
    (\kappa_1, 2k_2t)=1, \\
    \rad(s_1)\mid \kappa_1}}
    \! 2^{\omega(\kappa_1)}\mu(\kappa_1)^2
\eta(\kappa_1) = 
\frac{1}{2\pi i} \int_{c-iT}^{c+iT} \!\!F(s) \frac{y^s}{s} \mathrm{d}s + E_{\varepsilon,1},
\end{align*}
with $\max\{3,\rad(s_1)\}\le y\le Y/(k_2s_1)^2$, $c=1+\log(2y)^{-1}$ and
$
    E_{\varepsilon, 1} \ll_\varepsilon 2^{\omega(s_1)}
    (y^\varepsilon + \frac{y^{1+\varepsilon}}{T}).
$

Let $\eta^\star$ denote the primitive character
inducing $\eta$, with conductor $q_{\eta^\star}$.
Since the Euler factors
relating $L(s,\eta)$ and $L(s,\eta^\star)$
are of order $\ll_\varepsilon (XYZ)^\varepsilon$
uniformly for $\Re(s)\ge 1/2+\varepsilon$
by the restriction $|r|\le (XYZ)^2$,  the convexity
bound and \eqref{eq:facorisationF} imply \[
|L(s,\eta)| \ll_\varepsilon
(XYZ)^\varepsilon|L(s,\eta^\star)| \ll_\varepsilon (XYZ)^\varepsilon (q_{\eta^\star}T)^{
\max\{0,(1-\sigma)/2\}+\varepsilon}.
\]
On the
horizontal segments connecting $c\pm iT$
to $1/2+\varepsilon\pm iT$, this gives
\[
F(\sigma\pm iT)
\ll_\varepsilon
\frac{2^{\omega(s_1)}}{\rad(s_1)^\sigma}
(s_1k_2tq_{\eta^\star}Ty)^\varepsilon (XYZ)^\varepsilon  
(q_{\eta^\star}T)^{\max\{0,1-\sigma\}},
\]
uniformly for $1/2+\varepsilon\le\sigma\le c$.
Shifting the contour to the line $\Re(s) = 1/2+ \varepsilon$ therefore gives
\[
\sum_{\substack{1<\kappa_1 \le y, \\
    (\kappa_1, 2k_2t)=1, \\
    \rad(s_1)\mid \kappa_1}}
    \! 2^{\omega(\kappa_1)}\mu(\kappa_1)^2
\eta(\kappa_1) = 
    \frac{1}{2\pi i} \int_{\frac{1}{2}+\varepsilon-iT}^{\frac{1}{2}+\varepsilon+iT} \!\!F(s) \frac{y^s}{s} \mathrm{d}s +E_{\varepsilon,1} +E_{\varepsilon,2},
\]
where 
\[
    E_{\varepsilon,2} \ll 2^{\omega(s_1) } (s_1 k_2 t q_{\eta^\star} T XYZ)^\varepsilon
    \left(
    \frac{y}{\rad (s_1) T} + \frac{(y q_{\eta^\star})^{1/2}}{(\rad (s_1) T)^{1/2}}
    \right).
\]

By \eqref{eq:facorisationF}, we have 
\begin{align*}
    |F(\tfrac{1}{2}+\varepsilon+iv)|
\ll_\varepsilon
\frac{2^{\omega(s_1)}}{\rad(s_1)^{1/2}}
(s_1k_2tq_{\eta^\star}XYZ(1+|v|))^\varepsilon
|L(\tfrac{1}{2}+\varepsilon+iv,\eta^\star)|^2.
\end{align*}
Therefore, we may bound the integral 
\[
       \int_{\frac{1}{2}+\varepsilon-iT}^{\frac{1}{2}+\varepsilon+iT} F(s) \frac{y^s}{s} \mathrm{d}s \ll_\varepsilon 
       \frac{2^{\omega(s_1) } y^{1/2+\varepsilon}}{\rad(s_1)^{1/2} }
(s_1k_2tq_{\eta^\star} T XYZ)^\varepsilon
       \int_{-T}^{T} \frac{|L(\tfrac{1}{2}+\varepsilon+iv,\eta^\star)|^2}{1+|v|} \mathrm{d}v .
\]
Choosing $T=2yq_{\eta^\star}\rad(s_1)$, we get
\[
E_{\varepsilon,2}
\ll_\varepsilon
2^{\omega(s_1)}
(s_1k_2tq_{\eta^\star} TXYZ)^\varepsilon
\left(
\frac{1}{q_{\eta^\star}\rad(s_1)^2}
+
\frac{1}{\rad(s_1)}
\right).
\]
Since $T\le 2Yq_{\eta^\star}$, setting $T_\eta=2Yq_{\eta^\star}$, the whole contribution can thus be bounded uniformly in $y$ by 
\begin{align*}
\sum_{\substack{1<\kappa_1\le y\\
(\kappa_1,2k_2t)=1\\
\rad(s_1)\mid\kappa_1}} \!\!\!\!
2^{\omega(\kappa_1)}\mu(\kappa_1)^2\eta(\kappa_1)
\ll_\varepsilon
\frac{2^{\omega(s_1)}y^{1/2+\varepsilon}}
{\rad(s_1)^{1/2}}
(s_1k_2tq_{\eta^\star} T_{\eta} XYZ)^\varepsilon
\bigg( \!
1 \!
+\!
\int_{-T_\eta}^{T_\eta}
\frac{|L(\frac{1}{2}+\varepsilon+iv,\eta^\star)|^2}{1+|v|}
\mathrm{d}v
\bigg).
\end{align*}
Thus, by partial summation, the sum in \eqref{eq:kappa1sumtobound}
can be bounded as 
\[
\ll_\varepsilon
\frac{2^{\omega(s_1)}}{\rad(s_1)^{1/2}}
\bigl(s_1k_2tq_{\eta^\star}T_\eta XYZ\bigr)^\varepsilon
\frac{Z\tau\kappa_2t^2}{X|r|}
\left(
1+
\int_{-T_\eta}^{T_\eta}
\frac{|L(\tfrac12+\varepsilon+iv,\eta^\star)|^2}{1+|v|}
\mathrm{d}v
\right).
\]
Here we used the bounds $ 
\left|\widehat{w} (v) \right|\ll \min(1,|v|^{-1})
$ and $\left|\widehat{w}'(v)\right|
\ll \min(1,Z|v|^{-2})$ from \eqref{eq:fouriertrafobound} and \eqref{eq:whatderivativebound}.

For the sum over $s_1$, we observe that 
\[
\sum_{\substack{s_1\ge1\\2\nmid s_1}}
\frac{2^{\omega(s_1)}}{\rad(s_1)^{1/2}}
\prod_{p\mid s_1}
\frac{1+\nu_p(s_1)}{p^{2\nu_p(s_1)}}
=
\prod_{p>2}
\left(
1+\frac{2}{p^{1/2}}
\sum_{a\ge1}\frac{a+1}{p^{2a}}
\right)
\ll 1.
\]
Thus it remains to bound 
\begin{align*}
Z(XYZ)^\varepsilon
\sum_{k_2\le \sqrt Y} 
\frac{d(k_2^2)}{k_2^2}
\sum_{t\le \sqrt X}
\frac{2^{\omega(t)}\kappa_2\tau}{t}
\sum_{x_0\mid\tau}
x_0^{3/2} &\!\!\!
\sum_{\substack{r\neq 0, \\ |r|\le (XYZ)^2 \\ q\ne \kappa_2\tau/x_0}}\!\! 
\frac{|C_{k_2}(r)|}{|r|} \prod_{p\mid t}\frac1{p(p^2-1)}
\prod_{\substack{p\mid t,\\p\mid r}}(p^3-1) \\
& \qquad \cdot \left( 
1+ 
\int_{-T_\eta}^{T_\eta}
\frac{|L(\tfrac12+\varepsilon+iv,\eta^\star)|^2}{1+|v|}
\mathrm{d}v
\right),
\end{align*}
where, as always, $r=\sigma 2^j q m^2.$

The conductor of $\eta^\star$ is of order $d:= (q \kappa_2 \tau/x_0)_{\square \mathrm{ free}}$ (more concretely, we have $q_{\eta^\star}\in\{d,4d,8d\}$), where 
 $d$ denotes the squarefree part of $q \kappa_2 \tau/x_0$.
For any $d > 1$, we now bound from above the sum
\[
\sum_{\substack{k_2\le \sqrt Y,\ t\le \sqrt X,\ x_0\mid\tau, \, (k_2,t)=1\\
0<|r|\le (XYZ)^2, \,
(r,x_0)=1,\\
\\
d=(q\kappa_2\tau/x_0)_{\square\mathrm{free}}}} \!\!\!
\frac{d(k_2^2)}{k_2^2}
\frac{2^{\omega(t)}\kappa_2\tau}{t}
x_0^{3/2}
\frac{|C_{k_2}(r)|}{|r|}
\prod_{p\mid t}\frac1{p(p^2-1)}
\prod_{\substack{p\mid t,\\p\mid r}}(p^3-1).
\]
Since all summands are nonnegative, this is 
\[
\ll \sum_{\substack{\text{$p\mid k_2 t r$ implies $p \le (XYZ)^2$} \\
2\nmid k_2t,\ r\neq 0 \\
x_0 \mid \tau \\
(k_2,t)=(r,x_0)=1 \\
d=(q\kappa_2\tau/x_0)_{\square\mathrm{free}}}}\!\!\!
\frac{d(k_2^2)}{k_2^2}
\frac{2^{\omega(t)}\kappa_2\tau}{t}
x_0^{3/2}
\frac{|C_{k_2}(r)|}{|r|}
\prod_{p\mid t}\frac{1}{p(p^2-1)}
\prod_{\substack{p\mid t,\\p\mid r}}(p^3-1)
= 4 \prod_{\substack{2<p \le (XYZ)^2}} h_p(d),
\]
where the final expression follows by multiplicativity
for some function $h_p(d)$.

If $ p \nmid d$, $p \neq 2$, considering the four types of summands $p \nmid k_2t$, $ p \mid k_2, p \mid t \text{ but } p \nmid x_0$ and $ p \mid x_0$, we get 
\[
    h_p(d) = \sum_{m\ge 0} \frac{1}{p^{2m}} \Bigg( 1 + \sum_{j\ge 1} \frac{(2j+1)(p-1)}{p^{2j+1}} 
    + \sum_{i \ge 1}\frac{2(p^3-1)}{p^{i+1} (p^2 -1)}
    \Bigg) + \sum_{k\ge 1} \frac{2 p^{3/2}}{p^k(p^2-1)} =  1 + O(p^{-1}).
\]
If $ p \mid d$, we similarly get
\[
    h_p(d)  =\frac{1}{p(1-p^{-2})}
+\left(1+\frac{p-1}{p^2(1-p^{-2})}\right)
\sum_{j\ge1}\frac{2j+1}{p^{2j}}
+\frac{2}{(p-1)(p^2-1)}
\left(1+\frac{p^3-1}{p^2(1-p^{-2})}\right),
\]
and so $h_p(d) = O(p^{-1})$ in this case.
It follows that
\[
    \prod_{\substack{2<p \le (XYZ)^2}}
    h_p(d) \ll_\varepsilon (XYZ)^\varepsilon d^{-1}.
\]
For each fixed $d>1$, there are at most four possible
primitive characters $\eta^\star$ whose fundamental
discriminants have odd squarefree part $d$.

Splitting the integral into $[-1,1]$ and dyadic ranges $\pm [2^j,2^{j+1}]$ for $j\ll \log T$, we have
\begin{align*}
    \int_{-T}^{T}  \frac{|L(\tfrac{1}{2}+\varepsilon+iv,\eta^\star)|^2}{1+|v|} \mathrm{d}v
    &\ll  \int_{-1}^{1} {|L(\tfrac{1}{2}+\varepsilon+iv,\eta^\star)|^2} \mathrm{d}v 
    \\
    &+ \sum_{j \ll \log T}  \int_{\pm [2^j,2^{j+1}]} \frac{|L(\tfrac{1}{2}+\varepsilon+iv,\eta^\star)|^2}{2^j} \mathrm{d}v .
\end{align*}
Summing over $\eta^\star$ with $R \le q_{\eta^\star} \le 2R$ (so that $T_\eta \ll RY$) yields by Lemma~\ref{lem:lfunctionmeansquare}
\[
    \sum_{ \eta^\star:\,R\le q_{\eta^\star} \le 2R }   (XYZ)^\varepsilon R^{-1+\varepsilon} \int_{-T_\eta}^{T_\eta}  \frac{|L(\tfrac{1}{2}+\varepsilon+iv,\eta^\star)|^2}{1+|v|} \mathrm{d}v \ll_\varepsilon (XYZ)^\varepsilon Y^{3 \varepsilon/2} R^{3\varepsilon}.
\]
Summing over all dyadic blocks with $R\ll (XYZ)^2\sqrt{XY}$ completes the proof.
\end{proof}

\bibliographystyle{amsplain}
\bibliography{refs.bib}

\end{document}